\documentclass[11pt]{article}
\usepackage[utf8]{inputenc}
\usepackage{graphicx}
\usepackage{caption}
\usepackage{amsmath}
\usepackage{amssymb}
\usepackage{float}
\usepackage{hyperref}
\usepackage{amsthm}
\usepackage{xcolor}
\usepackage{comment}
\usepackage{amsfonts}
\usepackage{dsfont}
\usepackage{stmaryrd}
\usepackage{xcolor}
\usepackage{subcaption}

\usepackage{comment} \usepackage{siunitx}
\usepackage{mathtools} 
\mathtoolsset{showonlyrefs} 

\usepackage{lineno} 

\def \Sum{\displaystyle\sum}

\def\Zc{{\cal Z}}

\def\Fc{{\cal F}}
\def\Uc{{\cal U}}
\def\Tc{{\cal T}}

\def\Lc{{\cal L}}
\def\Sc{{\cal S}}
\def\Xc{{\cal X}}

\def\Jc{{\cal J}}
\def\eps{{\varepsilon}}

\def\bZc{\boldsymbol{\Zc}}

\def\btheta{\boldsymbol{\theta}}

\newcommand{\R}{\mathbb{R}}
\newcommand{\F}{\mathbb{F}}

\newcommand{\N}{\mathbb{N}}

\def \S{\mathbb{S}}
\def \P{\mathbb{P}}
\newcommand{\E}{\mathbb{E}}

\newcommand{\di}{\mathrm{d}}
\def\Lc{{\cal L}}

\def \mrx{\mathrm{x}} 
\def \trans{^{\scriptscriptstyle{\intercal}}}
\def \ep{\hbox{ }\hfill$\Box$}

\DeclareMathOperator{\Tr}{Tr}

\def\beqs{\begin{eqnarray*}}
\def\enqs{\end{eqnarray*}}
\def\beq{\begin{eqnarray}}
\def\enq{\end{eqnarray}}

\usepackage[algo2e,ruled,vlined]{algorithm2e} 
 \usepackage{algorithm}
 \usepackage{algorithmicx}
\usepackage{algpseudocode}
\usepackage{ marvosym }
\usepackage{csquotes}
\usepackage{hyperref}

\def \eps{\varepsilon}

\def\argmin{\mathop{\rm argmin}}
\def\argmin_#1{\underset{#1}{\mathrm{argmin\, }}}

\def \ep{\hbox{ }\hfill$\Box$}

\newtheorem{Theorem}{Theorem}[section]

\newtheorem{Remark}[Theorem]{Remark}

\numberwithin{equation}{section} 

\usepackage[backend=bibtex,maxbibnames=5, style=alphabetic]{biblatex}

\def \trans{^{\scriptscriptstyle{\intercal}}}

\begin{document}
\title{Neural networks-based algorithms for stochastic control and PDEs in finance
\footnote{This paper is a contribution for the {\it Machine Learning for Financial Markets: a guide to contemporary practices}, 
Cambridge University Press, Editors: Agostino Capponi and Charles-Albert Lehalle. 
This study was supported by FiME (Finance for Energy Market Research Centre) and
the ``Finance et D\'eveloppement Durable - Approches Quantitatives'' EDF - CACIB Chair.
}
}

\author{Maximilien \textsc{Germain}
\footnote{EDF R\&D, LPSM, Université de Paris  \sf \href{mailto:Maximilien.Germain at edf.fr}{mgermain at lpsm.paris}} \and  Huyên \textsc{Pham}
\footnote{LPSM, Université de Paris, FiME, CREST ENSAE \sf \href{mailto:pham at lpsm.paris}{pham at lpsm.paris}} \and  Xavier \textsc{Warin}
\footnote{EDF R\&D, FiME \sf \href{mailto:xavier.warin at  edf.fr}{xavier.warin at edf.fr}} 
}
\maketitle              
\begin{abstract}
    This paper presents machine learning techniques and deep reinforcement lear\-ning-based algorithms for the efficient resolution of nonlinear partial differential equations and dynamic optimization problems  arising in  investment decisions and derivative pricing in financial engineering. We survey recent results in the literature, 
present new developments, notably in the fully nonlinear case, and compare the different schemes illustrated by numerical tests on various financial applications. We conclude by highlighting some future research directions. 

\end{abstract}

\section{Breakthrough in the  resolution of high dimensional non-linear problems}



The numerical resolution of control problems and nonlinear PDEs-- arising in several financial applications such as portfolio selection, hedging, or derivatives pricing--is subject to the so-called ``curse of dimensionality", making impractical the discretization of the state space in dimension greater than 3  by using classical PDE resolution methods such as finite differences schemes. 
Probabilistic regression Monte-Carlo methods based on a Backward Stochastic Differential Equation (BSDE) representation of semilinear PDEs have been developed in \cite{Z04}, \cite{BT04}, \cite{GLW05} to overcome this obstacle. These mesh-free techniques are successfully applied upon dimension 6 or 7, nevertheless, their use of regression methods requires a number of basis functions growing fastly  with the dimension. What can be done to further increase the dimension of numerically solvable problems?



A breakthrough with deep learning based-algorithms has been made in the last five years towards this computational challenge, and we mention the recent survey by \cite{bechutjenkuc20}. 
The main interest in the use of machine learning techniques for control and PDEs
is the ability of deep neural networks to efficiently represent high dimensional functions without using spatial grids, and with no curse of dimensionality \cite{grohorjenvon18}, \cite{hutetal20}. 
Although the use of neural networks for solving PDEs is not new,  see e.g. \cite{dispha94},  the approach has been successfully revived with new ideas and directions.  
Neural networks have known a increasing popularity since the works on Reinforcement Learning for solving the game of Go  by Google DeepMind teams. 
These empirical successes and the introduced methods allow to solve control problems in moderate or large dimension. Moreover, recently developed open source libraries like Tensorflow 
and  Pytorch 
also offer an accessible framework to implement these algorithms. 


 A first natural use of neural networks for stochastic control concerns the discrete time setting, with the study of Markov Decision Processes, either in a brute force fashion or by using dynamic programming approaches. 
 In the continuous time setting, and in  the context of PDE resolution, we present various methods. A first kind of schemes is rather generic and can be applied to a variety of PDEs coming from a large range of applications. Other schemes rely on BSDE representations, strongly linked to stochastic control problems. In both cases, numerical evidence seems to indicate that the methods can be used in large dimension, greater than 10 and up to 1000 in certain studies. Some theoretical results also illustrates the convergence of specific algorithms. These advances pave the way for new methods dedicated to the study of large population games, studied in the context of mean field games and mean field control problems. 

 The outline of this article is the following. We first focus on some schemes for discrete time control in Section \ref{sec: discrete time} before presenting generic machine learning schemes for PDEs in Subsection \ref{sec : deterministic}. Then we review BSDE-based machine learning methods for semilinear equations in Subsection \ref{sec: bsde semi schemes}. Existing algorithms for fully non-linear PDEs are detailed in Subsection \ref{sec: bsde fully schemes} before presenting new BSDE schemes designed to treat this more difficult case. Numerical tests on CVA pricing and portfolio selection are conducted in Section \ref{sec: numerics} to compare the different approaches. 
Finally, we highlight  in Section \ref{sec: extensions} further directions and perspectives including  recent advances for the resolution of mean field games and mean field control problems with or without model.

\section{Deep learning approach for stochastic control}\label{sec: discrete time}

We present in this section some recent breakthrough in the numerical resolution of stochastic control in high dimension by means of machine learning techniques. 
We consider a model-based setting in discrete-time, i.e., a Markov decision process, that could possibly be obtained from the time discretization of a continuous-time stochastic control problem. 
 
Let us fix  a probability space $(\Omega,\Fc,\P)$ equipped with a filtration $\F$ $=$ $(\Fc_t)_t$ representing the available information at any time $t$ $\in$ $\N$ ($\Fc_0$ is  the trivial $\sigma$-algebra).  
The evolution of the system is described by a model dynamics for the state process $(X_t)_{t\in \N}$ valued in $\Xc$ $\subset$ $\R^d$: 
\begin{align} \label{dynXdis} 
X_{t+1} &= \;  F(X_t,\alpha_t,\eps_{t+1}), \quad t \in \N, 
\end{align} 
where $(\eps_t)_t$ is a sequence of i.i.d. random variables valued in $E$, with $\eps_{t+1}$  $\Fc_{t+1}$-measurable containing all the noisy information arriving between $t$ and $t+1$,  and 
$\alpha$ $=$ $(\alpha_t)_t$ is the control  process valued in $A$ $\subset$ $\R^q$.  The dynamics function $F$ is a measurable function from $\R^d\times\R^q\times E$ into $\R^d$, and assumed to be known.  
Given a running cost function $f$,  
a finite horizon $T$ $\in$ $\N^*$,  and a terminal cost function, 
the problem is to minimize over control process $\alpha$ a functional cost
\begin{align} \label{Jdis} 
J(\alpha) &= \; \E \Big[ \sum_{t=0}^{T-1} f(X_t,\alpha_t) + g(X_T) \Big]. 
\end{align}
In some relevant  applications, we may require constraints on the state and control in the form: $(X_t,\alpha_t)$ $\in$ $\Sc$, $t\in\N$. 
for some subset $\Sc$ of $\R^d\times\R^q$. This can be handled by relaxing the state/constraint and introducing into the costs a penalty function $L(x,a)$: $f(x,a)$ $\leftarrow$ $f(x,a) + L(x,a)$, and $g(x)$ $\leftarrow$ $g(x) + L(x,a)$.  For example, if the constraint set is in the form: 
$\Sc$ $=$ $\{(x,a) \in \R^d\times\R^q: h_k(x,a) = 0, k=1,\ldots,p, h_k(x,a) \geq 0, k=p+1,\ldots,m\}$, then one can take as penalty functions: 
\begin{align}
L(x,a) &= \; \sum_{k=1}^p \mu_k |h_k(x,a)|^2+ \sum_{k=p+1}^m \mu_k \max(0,-h_k(x,a)), 
\end{align}
where $\mu_k$ are penalization parameters (large in practice) see e.g. \cite{hanE}.

\subsection{Global approach}

The method consists simply in  approximating at any time $t$, the feedback control, i.e. a function of the state process, by a neural network (NN):  
\begin{align}
\alpha_t & \simeq \;  \pi^{\theta_t}(X_t), \quad t = 0,\ldots,T-1,  
\end{align}
where $\pi^\theta$ is a feedforward neural network on $\R^d$  with parameters $\theta$, and then to minimize over the global set of parameters $\btheta$ $=$ $(\theta_0,\dots,\theta_{T-1})$ the quantity (playing the role of loss function)
\begin{eqnarray*} 
\tilde J(\btheta)  & = &  \E \Big[  \sum_{t=0}^{T-1}  f(X_t^{\btheta},\pi^{\theta_t}(X_t^{\btheta}) )   + g(X_T^{\btheta}) \Big], 
\end{eqnarray*} 
where $X^{\btheta}$ is the state  process associated with the NN feedback controls:
\begin{eqnarray*}
X_{t+1}^{\btheta} &= & F(X_t^{\btheta},\pi^{\theta_t}(X_t^{\btheta}), \eps_{t+1}), \quad t = 0, \ldots,T-1. 
\end{eqnarray*}
This basic idea of approximating control by parametric function of the state was  proposed in \cite{gobmun05}, and updated  with the use of (deep) neural networks by  \cite{hanE}.  
This method met  success due to its simplicity and the easy accessibility of  common libraries like TensorFlow for optimizing the parameters of the neural networks.  
Some recent extensions of this approach dealt with stochastic control problems with delay, see \cite{hanhu21}. 
However, such global optimization over a huge set of parameters $\btheta$ $=$ $(\theta_0,\dots,\theta_{T-1})$ may suffer from being stuck in suboptimal traps and thus does not converge, especially  for large horizon $T$.  An alternative is to consider controls $\alpha_t  \simeq \;  \pi^{\theta}(t,X_t),$ $ t = 0,\ldots,T-1, $ with a single neural network $\pi^{\theta}$ giving  more stabilized results as studied by \cite{FMW19}. We focus here on feedback controls, which is not  a restriction as we are in a Markov setting. For path-dependent control problems, we may consider  recurrent neural networks to take into consideration the past of state trajectory as input of the policy.

\subsection{Backward dynamic programming approach}

In  \cite{bachurlanpha19b}, the authors propose methods  that combine ideas from numerical probability and  deep reinforcement learning. Their algorithms are based on the classical dynamic programming (DP), (deep) neural networks for  the approximation/learning of the optimal policy and value function, 
and Monte-Carlo regressions with performance and value iterations.

The first algorithm, called NNContPI, is a combination of dynamic pro\-gramming  and the approach in \cite{hanE}. It learns sequentially the control by NN $\pi^\theta(.)$ and performance iterations, and is designed as follows: \\
\begin{algorithm2e}[H]
 	\caption{NNContPI}
 	\label{algo:nncontPI}
 	\textbf{Input:} the training distributions $(\mu_t)_{t=0}^{T-1}$\;
 	\textbf{Output:} estimates of the optimal strategy $(\hat\pi_t)_{t=0}^{T-1}$\;
 	\For{$t$ $=$ $T-1,\ldots,0$}{
 		Compute $\hat\pi_t := \pi^{\hat\theta_t}$ with 
 \begin{align}
\hat \theta_t  & \in  \;    {\rm arg}\min_{\theta} \E \Big[ f\big(X_t,\pi^\theta(X_t)\big) +  \sum_{s=t+1}^{T-1} 
 f\big(X^{\theta}_s, \hat{\pi}_s \big( X^{\theta}_s \big)\big)  + g\big(  X^{\theta}_T  \big)\Big]  \\
 \label{computa}
 \end{align}
where $X_t \sim \mu_t$ and where $\big(X^{\theta}_s\big)_{s=t+1}^T$ is defined by induction as:
 \begin{equation*}
 \left\{
 \begin{array}{ccl}
 X^{\theta}_{t+1} &=& F\big( X_t, \pi^\theta(X_t), \eps_{t+1} \big),  \\
 X^{\theta}_{s+1} &=& F\big( X^{\theta}_{s}, \hat \pi_s(X^{\theta}_{s}), \eps_{s+1} \big), \;\;\; \mbox{  for } s= t+1, \ldots, T-1.
 \end{array}
 \right.
 \end{equation*}
 }
\end{algorithm2e}

The second algorithm, refereed to as Hybrid-Now,   combines optimal policy estimation by neural networks and dynamic programming principle, and relies on  an hybrid procedure between value and performance iteration to approximate the value function by 
neural network $\Phi^\eta(.)$ on $\R^d$  with parameters $\eta$. \\
 \begin{algorithm2e}[H]
 	\caption{Hybrid-Now}
 	\label{algo:Hybrid-Now}
 	\textbf{Input:} the training distributions $(\mu_t)_{t=0}^{T-1}$\;
 	\textbf{Output:} \\
 	-- estimate of the optimal strategy $(\hat\pi_t)_{t=0}^{T-1}$\;
 	-- estimate of the value function $(\hat V_t)_{t=0}^{T-1}$\;
 Set $\hat V_T$ $=$ $g$\;
 	\For{$t$ $=$ $T-1,\ldots,0$}{
 Compute: 
\begin {align}
\hat \theta_t  & \in  \;     {\rm arg}\min_{\theta} \E \Big[ f\big(X_t^{},\pi^\theta(X_t)\big)+ \hat{V}_{t+1}(X_{t+1}^{\theta}) \Big] \label{computanow}
\end{align}
where $X_t \sim \mu_t$, and  $X_{t+1}^{\theta}$ $=$ $F(X_t^{},\pi^\theta(X_t^{}),\eps_{t+1}^{})$\;
 		Set $ \hat \pi_t  :=  \pi^{\hat\theta_t}$; \Comment{$\hat \pi_t$ is the estimate of the optimal policy at time $t$}\\
 Compute
 \begin{align}
 \hat\eta_t & \in  {\rm arg}\min_{\eta} \E  \Big|  f(X_t^{},\hat \pi_t(X_t^{}))+ \hat{V}_{t+1}(X_{t+1}^{\hat\theta_t}) - \Phi^\eta(X_t^{}) \Big|^2.  \label{computv}
 \end{align} 
 Set $\hat V_t = \Phi^{\hat\eta_t}$; \Comment{$\hat V_t$ is the estimate of the value function at time $t$}
 	}
 \end{algorithm2e}

The convergence analysis  of  Algorithms  NNContPI and Hybrid-Now  are stu\-died in \cite{huretal18}, and various applications in finance are implemented in  \cite{bachurlanpha19b}. 
These algorithms  are well-designed for control problems with continuous control space $A$ $=$ $\R^q$ or a ball in $\R^q$. In the case where the control space $A$ is finite, it is relevant to randomize controls, and then use classification methods by approximating the distribution of controls with neural networks and  Softmax activation function.

\section{Machine learning algorithms  for nonlinear PDEs}

By change of time scale, Markov decision process   \eqref{dynXdis}-\eqref{Jdis} can be obtained from  the time discretization of a continuous-time stochastic control problem with controlled diffusion dynamics on $\R^d$
\begin{align}
dX_t &=\;  b(X_t,\alpha_t) dt + \sigma(X_t,\alpha_t) dW_t,    
\end{align}
and cost functional to be minimized over control process $\alpha$ valued in $A$
\begin{align}
J(\alpha) &= \;  \E\Big[ \int_0^T f(X_t,\alpha_t) dt + g(X_T) \Big].    
\end{align}
In this case, it is well-known, see e.g. \cite{pha09}, that the dynamic programming Bellman equation leads to a  PDE  in the form
\begin{equation}
\begin{cases}
\partial_t u + H(x,D_x u,D_x^2 u)  \; =\;  0, \mbox{ on } [0,T)\times\R^d \\
u(T,.) \; =\; g \quad \mbox{ on } \R^d,
\end{cases}
\end{equation}
where $H(x,z,\gamma)$ $=$ $\inf_{a \in A}\big[b(x,a).z + \frac{1}{2}{\rm tr}(\sigma\sigma\trans(x,a)\gamma)  + f(x,a)\big]$ is the so-called Hamiltonian function.  The numerical resolution of such class of second-order parabolic PDEs  will be addressed in this section.

\subsection{Deterministic approach by neural networks}\label{sec : deterministic}

In the schemes below, differential operators are evaluated by automatic differentiation of the network function approximating the solution of the PDE. Machine learning libraries such as Tensorflow or Pytorch allow to efficiently compute these derivatives. The studied PDE problem is 
\begin{equation}\label{eq: galerkin}
    \begin{cases}
        \partial_t u + \Fc u \; = \;  0 &\mathrm{ on } \ [0,T)\times\Lambda\\
        u(T,\cdot) \; = \;  g  &\mathrm{ on } \ \Lambda\\
        u(t,x) \; = \;  h(t,x)  &\mathrm{ on } \ [0,T)\times\partial\Lambda,
    \end{cases}
\end{equation}
with $\Fc$ a space differential operator, $\Lambda$ a subset of $\R^d$.

 \vspace{1mm}
 
\noindent $\bullet$ {\bf Deep Galerkin Method \cite{SS17}.}

\vspace{1mm}
 
\noindent The Deep Galerkin Method is a meshfree machine learning algorithm to solve PDEs on a domain, eventually with boundary conditions. The principle is to sample time and space points according to a training measure, e.g. uniform on a bounded domain,
and minimize a performance measure quantifying how well a neural network satisfies the differential operator and boundary conditions. The method consists in minimizing over neural network $\Uc$ $:$ $\R\times\R^d$ $\rightarrow$ $\R^d$,  the $L^2$ loss
\begin{equation}
    \E|\partial_t \Uc(\tau,\kappa) + \Fc \Uc(\tau,\kappa)|^2 +  \E|\Uc(T,\xi)-g(\xi)|^2 +  \E|\Uc(\tau,\kappa) - h(\tau,\kappa)|^2
\end{equation} with  $\kappa,\tau,\xi$ independent random variables in $\Lambda \times [0,T) \times \partial \Lambda$.  \cite{SS17} also prove a convergence result (without rate) for the Deep Galerkin method. This method is tested on financial problems by \cite{alaradi2018solving}. A major advantage to this method is its adaptability to a large range of PDEs with or without boundary conditions. Indeed the loss function is straightforwardly modified according to changes in the constraints one wishes to enforce on the PDE solution. A related approach is the deep parametric PDE method, see \cite{kholuyin20}, and \cite{glawun20} applied to option pricing.
\vspace{2mm}

\noindent $\bullet$ {\bf Other approximation methods}

\begin{itemize}
\item[(i)] {\it Physics informed neural networks} \cite{raissiComp}. Physics informed neural networks use both data (obtained for a limited amount of samples from a PDE solution), and theoretical dynamics to reconstruct solutions from PDEs. The convergence of this method in the Second Order linear parabolic (or elliptic) case is proven in \cite{SDK20}, see also \cite{GAS20}.
\item[(ii)] {\it Deep Ritz method}  \cite{ritz}.  The Deep Ritz method focuses on the resolution of the variational formulation from elliptic problems where the integral is evaluated by randomly sampling time and space points, like in the Deep Galerkin method \cite{SS17} and the minimization is performed over the parameters  of a neural network.  This scheme is tested on Poisson equation with different types of boundary conditions. \cite{MZ19} studies the convergence of the Deep Ritz algorithm.
\end{itemize}

\subsection{Probabilistic approach by neural networks}

\subsubsection{Semi-linear case}\label{sec: bsde semi schemes}

In this paragraph, we consider semilinear PDEs of the form \begin{equation}\label{eq: semilinear PDE}
    \begin{cases}
        \partial_t u + \mu \cdot D_x u + \frac{1}{2} \Tr(\sigma\sigma\trans D^2_x u) \; = \;  f(\cdot,\cdot,u,\sigma\trans D_x u)  &\mathrm{ on } \ [0,T)\times\R^d\\
        u(T,\cdot) \; = \;  g  &\mathrm{ on } \ \R^d.
    \end{cases}
\end{equation}
for which we have the  forward backward SDE representation 
\begin{equation} \label{BSDEfeyn} 
\left\{
\begin{array}{ccl}
 Y_t &=& g(\Xc_T) - \int_t^T f(s,\Xc_s,Y_s,Z_s) \di s - \int_t^T  Z_s \cdot  \di W_s, \quad 0 \leq t \leq T, \\
 \Xc_t &=& \Xc_0 + \int_0^t \mu(s,\Xc_s) \di s + \int_0^t \sigma(s,\Xc_s) \di W_s, 
 \end{array}
 \right.
\end{equation} 
via the (non-linear) Feynman-Kac formula: $Y_t$ $=$ $v(t,X_t)$, $Z_t$ $=$ $\sigma\trans(t,\Xc_t)D_x v(t,X_t)$, $0\leq t\leq T$, see \cite{PP90}.  

Let $\pi$ be a subdivision $\{t_0=0<t_1<\cdots<t_N=T\}$ with modulus  $|\pi| := \sup_i \Delta t_i$, $\Delta t_i$ $:=$ $t_{i+1} - t_i$, 
satisfying $|\pi| = O\left(\frac{1}{N}\right)$, and consider the Euler-Maruyama 
discretization $(X_{i})_{i=0,\ldots,N}$ defined by
\begin{align} \label{eq: Euler}
  X_{i} & = \Xc_0 +  \sum_{j=0}^{i-1}  \mu(t_j,X_{j}) \Delta t_j +  \sum_{j=0}^{i-1}    \sigma(t_j,X_{j}) \Delta W_j,
\end{align} 
where  $\Delta W_j$ $:=$ $W_{t_{j+1}} - W_{t_j}$, $j$ $=$ $0,\ldots,N$. 
Sample paths of $(X_i)_i$ 
act as training data in the machine learning setting. Thus our training set can be chosen as large as desired, which is 
relevant  for training purposes as sit does not  lead to overfitting.

The time discretization of the BSDE \eqref{BSDEfeyn} can be written in backward induction as 
\begin{align} \label{Yinduc}
Y_{i}^\pi & = \; Y_{i+1}^\pi - f(t_i,X_i,Y_{i}^\pi,Z_{i}^\pi) \Delta t_i  - Z_{i}^\pi. \Delta W_{i}, \quad i=0,\ldots,N-1,  
\end{align}
which can be described as conditional expectation formulae
\begin{equation} \label{Ycond}
\left\{
\begin{array}{ccl}
Y_{i}^\pi & = &  \E_i \Big[ Y_{{i+1}}^\pi - f(t_i,X_i,Y_{i}^\pi,Z_{i}^\pi) \Delta t_i  \Big]  \\
Z_{i}^\pi & = &  \E_i \Big[  \frac{\Delta W_i}{\Delta t_i} Y_{{i+1}}^\pi \Big], \quad \quad i=0,\ldots,N-1, 
\end{array}
\right.
\end{equation}
where $\E_i$ is a notation for the conditional expectation w.r.t. $\mathcal{F}_{t_i}$.  

\vspace{3mm}

\noindent  $\bullet$ {\bf Deep BSDE scheme  \cite{Ehanjen17},  \cite{HJE17}.} 

\vspace{1mm}

\noindent The essence of this method is to write down the backward equation \eqref{Yinduc} as a forward equation. One appro\-ximates the initial condition $Y_0$ and the $Z$ component at each time  by networks functions taking the forward process
$X$ as input.  The objective function to optimize is the error between the reconstructed dynamics and the true terminal condition. More precisely, the problem is to minimize over  network functions $\Uc_0$ $:$ $\R^d$ $\rightarrow$ $\R$, and  sequences of network  functions 
$\bZc$ $=$ $(\Zc_i)_i$, $\Zc_i$ $:$ $\R^d$ $\rightarrow$  $\R^d$, $i$ $=$ $0,\ldots,N-1$,  the global quadratic loss function
\begin{align}
J_G(\Uc_0,\bZc) &=\; \E \Big| Y_N^{\Uc_0,\bZc} - g(X_N) \Big|^2,
\end{align}
where $(Y_i^{\Uc_0,\bZc})_i$ is defined by forward induction as 
\begin{eqnarray*}
Y_{i+1}^{\Uc_0,\bZc} & = & Y_{i}^{\Uc_0,\bZc} +  f(t_i,X_i,Y_i^{\Uc_0,\bZc},\Zc_i(X_i)) \Delta t_i + \Zc_i(X_i).\Delta W_i, \;  i = 0,\ldots,N-1, 
\end{eqnarray*}
starting from $Y_0^{\Uc_0,\bZc}$ $=$ $\Uc_0(\Xc_0)$. 
The output of this scheme, for the solution $(\widehat{\Uc}_0,\widehat\bZc)$ to this global minimization problem, supplies an approximation 
$\widehat{\Uc}_0$ of the solution $u(0,.)$ to the PDE at time $0$, and approximations $Y_i^{\widehat{\Uc}_0,\widehat\bZc}$  of the solution to the PDE \eqref{eq: semilinear PDE} at times $t_i$ evaluated at 
$\Xc_{t_i}$, i.e., of   $Y_{t_i}$ $=$ $u(t_i,\Xc_{t_i})$, $i$ $=$ $0,\ldots,N$. The convergence of this algorithm through a posteriori error is studied by \cite{HL18}, see also \cite{JL21}.
A variant is proposed by \cite{CWNMW19} which introduces a single neural network $\Zc(t,x) : [0,T]\times \R^d \mapsto \R^d$ instead of $N$ independent neural networks. This simplifies the optimization problem and leads to more stable solutions. A close method introduced by  \cite{raissi2018forwardbackward} uses also a single neural network $\Uc(t,x) : [0,T]\times \R^d \mapsto \R$ and estimates $Z$ as the automatic derivative in space of $\Uc$. 
We also refer to \cite{jacoum19} for a variation of this deep BSDE scheme to curve-dependent {PDE}s arising in the pricing under rough volatility model, to \cite{nusric20} for approximations methods for Hamilton-Jacobi-Bellman PDEs, and to \cite{kresteszo20} for extension of deep BSDE scheme to 
elliptic PDEs with applications in insurance. 

\vspace{3mm}

\noindent $\bullet$   \textbf{Deep Backward Dynamic Programming (DBDP) \cite{HPW19}.}

\vspace{1mm}

\noindent  The method builds upon the backward dynamic programming relation \eqref{Yinduc}  stemming from the time discretization of the BSDE, and approximates simultaneously at each time step $t_i$ the processes $(Y_{t_i},Z_{t_i})$ with neural networks trained with the forward diffusion process $X_i$ as input.  
The  scheme can be implemented in two similar versions: 
\begin{itemize}
\item[1.] {\it  DBDP1}. Starting from $\widehat{\Uc}_N^{(1)}$ $=$ $g$, proceed by backward induction for $i$ $=$ $N-1,\ldots,0$, by minimizing over network functions 
$\Uc_i$ $:$ $\R^d$ $\rightarrow$ $\R$, and $\Zc_i$ $:$ $\R^d$ $\rightarrow$ $\R^d$ the quadratic  loss function 
\begin{eqnarray*}
&  & J_i^{(B1)}(\Uc_i,\Zc_i)\\ 
&=& \E \Big| \widehat{\Uc}_{i+1}^{(1)}(X_{i+1}) -  \Uc_i(X_i) -  f(t_i,X_i,\Uc_i(X_i),\Zc_i(X_i)) \Delta t_i -  \Zc_i(X_i).\Delta W_i \Big|^2, 
\end{eqnarray*}
and update $(\widehat{\Uc}_i^{(1)},\widehat{\Zc}_i^{(1)})$ as the solution to this local minimization problem. 
\item[2.] {\it  DBDP2}. Starting from $\widehat{\Uc}_N^{(2)}$ $=$ $g$, proceed by backward induction for $i$ $=$ $N-1,\ldots,0$, by minimizing over $C^1$ network functions  
$\Uc_i$ $:$ $\R^d$ $\rightarrow$ $\R$ the quadratic loss function 
\begin{align}
& J_i^{(B2)}(\Uc_i) \\&=\; \E \Big| \widehat{\Uc}_{i+1}^{(2)}(X_{i+1}) -  \Uc_i(X_i) -  f(t_i,X_i,\Uc_i(X_i), \sigma(t_i,X_i)\trans D_x\Uc_i(X_i)) \Delta t_i \\
&  \quad \quad \quad - \;   D_x\Uc_i(X_i)\trans\sigma(t_i,X_i)\Delta W_i \Big|^2, 
\end{align}
where $D_x \Uc_i$ is the automatic differentiation of the network function $\Uc_i$. 
Update $\widehat{\Uc}_i^{(2)}$ as  the solution to this local minimization problem, and set 
$\widehat{\Zc}_i^{(2)}$ $=$ $\sigma\trans(t_i,.)D_x{\Uc}_i^{(2)}$. 
\end{itemize}
The output of DBDP  supplies an approximation $(\widehat{\Uc}_i^{},\widehat{\Zc}_i^{})$  
of the solution $u(t_i,.)$ and its gradient $\sigma\trans(t_i,.)D_xu(t_i,.)$  to the PDE \eqref{eq: semilinear PDE} on the time grid $t_i$, $i$ $=$ $0,\ldots,N-1$. 
The study of the approximation error due to the time discretization and the choice of the loss function is accomplished in \cite{HPW19}.

\vspace{3mm}


\noindent {\bf Variants and extensions of DBDP schemes}

\begin{itemize}
\item[(i)] A regression-based machine learning scheme inspired by regression Monte-Carlo methods 
for numerically computing condition expectations in the time discretization  \eqref{Ycond} of  the BSDE,  is given by: 
starting from $\hat\Uc_N$ $=$ $g$, proceed by backward induction for $i$ $=$ $N-1,\ldots,0$, in two regression  problems:  
\begin{itemize}
\item[(a)] Minimize over network functions $\Zc_i$ $:$ $\R^d$ $\rightarrow$ $\R^d$ 
\begin{align}
J_{i}^{r,Z}(\Zc_i) & = \; \E \Big| \frac{\Delta W_i}{\Delta t_i} \widehat{\Uc}_{i+1}(X_{i+1}) -  \Zc_i(X_i) \Big|^2 
\end{align}
and update $\widehat{\Zc}_i$ as  the solution to this minimization problem
\item[(b)]  Minimize over network functions $\Uc_i$ $:$ $\R^d$ $\rightarrow$ $\R$
\begin{align}
J_{i}^{r,Y}(\Uc_i) &= \; \E \Big|  \widehat{\Uc}_{i+1}^{}(X_{i+1}) -  \Uc_i(X_i) -  f(t_i,X_i,\Uc_i(X_i),\widehat{\Zc}_i(X_i)) \Delta t_i   \Big|^2
\end{align} 
and update $\widehat{\Uc}_i$ as  the solution to this minimization problem. 
\end{itemize}
Compared to these regression-based schemes, the DBDP scheme simultaneously estimates the pair component $(Y,Z)$ through the minimization of the loss functions $J_i^{(B1)}(\Uc_i,\Zc_i)$ (or 
$J_i^{(B2)}(\Uc_i)$ for the second version), $i$ $=$ $N-1,\ldots,0$. Interestingly, the convergence of the DBDP scheme can be confirmed by computing at each time step the infimum of loss function, which should vanish for the exact solution (up to the time discretization). In contrast, the infimum of the  loss functions in usual regression-based schemes is unknown for the true solution as it is supposed to match the residual of $L^2$-projection. Therefore the scheme accuracy cannot be directly verified. 
\item[(ii)] The DBDP  scheme is based on local resolution, and was first used to solve linear PDEs, see \cite{SSS18}. 
It is  also suitable to solve variational inequalities  and can be used to valuate American options as shown in \cite{HPW19}. 
Alternative methods consists in using the Deep Optimal Stopping scheme  \cite{becchejen19} or the method from \cite{becker2019stopping}. Some tests on Bermudan options are also performed by \cite{bermudan} and \cite{fujtaktak19} with some refinements of the  Deep BSDE scheme. 
\item[(iii)] The {\bf Deep Splitting (DS) scheme in \cite{BBCJN19} }  combines ideas from the DBDP2 and regression-based schemes. Indeed the current regression-approximation on $Z$ is estimated by the automatic differentiation of the neural network computed at the previous optimization step. The current approximation of $Y$ is then computed by a regression-type optimization problem. It can be seen as a local version of the global algorithm from \cite{raissi2018forwardbackward} or as a step by step Feynman-Kac approach. As the scheme is a local one, it can be used to valuate American options. 
The convergence of this method is studied by \cite{GPW20}.
\item[(iv)] Local resolution permits to add other constraints such as constraints on a replication portfolio using facelifting techniques as in \cite{kharroubi2020discretization}.
\item[(v)] The \textbf{Deep Backward Multistep (MDBDP) scheme \cite{GPW20}}  is described as follows: for $i$ $=$ $N-1,\ldots,0$, minimize 
over network functions $\Uc_i$ $:$ $\R^d$ $\rightarrow$ $\R$, and $\Zc_i$ $:$ $\R^d$ $\rightarrow$ $\R^d$ the loss function 
\begin{align} 
& J_i^{MB}(\Uc_i,\Zc_i) \\& = \; \E \Big| g(X_N) - \sum_{j=i+1}^{N-1} f(t_j,X_j,\widehat{\Uc}_j(X_j),\widehat{\Zc}_j(X_j))\Delta t_j - \sum_{j=i+1}^{N-1} \widehat{\Zc}_j(X_j). \Delta W_j  \nonumber \\
& \quad \quad \quad - \;  \Uc_i(X_i) -  f(t_i,X_i,\Uc_i(X_i),\Zc_i(X_i)) \Delta t_i -  \Zc_i(X_i).\Delta W_i \Big|^2  \label{JMB}
\end{align}
and update $(\widehat{\Uc}_i^{},\widehat{\Zc}_i^{})$ as the solution to this  minimization problem.  
This output  provides an approxi\-mation $(\widehat{\Uc}_i,\widehat{\Zc}_i)$ of the solution $u(t_i,.)$ to the PDE \eqref{eq: semilinear PDE} at times $t_i$, $i$ $=$ $0,\ldots,N-1$.

MDBDP  is a machine learning version of the Multi-step Forward Dynamic Programming method studied by \cite{BD07} and \cite{GT14}.  
Instead of solving at each time step two regression problems, our approach allows to consider only a single minimization as in the DBDP scheme. Compared to the latter, the multi-step consideration is expected to provide better accuracy by reducing the propagation of errors in the backward induction as it can be shown  comparing the error estimated in \cite{GPW20} and \cite{HPW19} both at theoretical and numerical level.
\end{itemize}

\subsubsection{Case of fully non-linear PDEs}\label{sec: bsde fully schemes}

In this paragraph, we consider fully non-linear PDEs in the form \begin{equation}\label{eq: fully nonlinear PDE}
    \begin{cases}
        \partial_t u + \mu \cdot D_x u + \frac{1}{2} \Tr(\sigma\sigma\trans D^2_x u) \; =  \;  F(\cdot,\cdot,u,D_x u, D_x^2 u)  &\mathrm{ on } \ [0,T)\times\R^d\\
        u(T,\cdot) \; = \;  g &\mathrm{ on } \ \R^d,
    \end{cases}
\end{equation} 
For this purpose, we introduce a forward diffusion process  $\Xc$ in $\R^d$ as in \eqref{BSDEfeyn}, and associated to the linear part $\Lc$  of the differential operator in the l.h.s. of  the PDE \eqref{eq: fully nonlinear PDE}.  
Since the function $F$ contains the dependence both on the gradient $D_x u$ and the Hessian $D_x^2 u$, we can shift 
the linear differential operator (left hand side) of the PDE \eqref{eq: fully nonlinear PDE} into the function $F$. However, in practice, this linear differential operator associated to a diffusion process $\Xc$ is used for training simulations in SGD of machine learning schemes. 
We refer to Section 3.1 in \cite{phawar19} for a discussion on the choice of the parameters $\mu$, $\sigma$. 
In the sequel, we assume for simplicity  that  $\mu$ $=$ $0$, and $\sigma$ is a constant invertible matrix.

Let us derive formally a BSDE representation for the nonlinear PDE \eqref{eq: fully nonlinear PDE}  on which we shall rely for designing our machine learning algorithm.
Assuming  that the solution $u$ to this PDE is smooth $C^2$, and denoting by $(Y,Z,\Gamma)$ the triple  of $\F$-adapted processes valued in 
$\R\times\R^d\times\S^d$,  defined by 
\begin{align}
Y_t \; = \; u(t,\Xc_t), \quad Z_t \; = \; D_x u(t,\Xc_t), \quad \Gamma_t \; = \; D_x^2 u(t,\Xc_t), \quad 0 \leq t \leq T, 
\end{align}
a direct application of It\^o's formula to $u(t,\Xc_t)$, yields that $(Y,Z,\Gamma)$ satisfies the backward equation
\begin{align} \label{BSDE2}
Y_t & = \; g(\Xc_T) - \int_t^T F(s,\Xc_s,Y_s,Z_s,\Gamma_s) \di s - \int_t^T    Z_s\trans \sigma  \di W_s, \; 0 \leq t \leq T.  
\end{align}



Compared to the case of semi-linear PDE of the form \eqref{eq: semilinear PDE}, the key point is the approximation/learning of the Hessian matrix $D_x^2 u$, hence of the $\Gamma$-component of the BSDE \eqref{BSDE2}.  
We present below different approaches for the approximation of the $\Gamma$-component. To the best of our knowledge, no theoretical convergence result is available for machine learning schemes in the fully nonlinear case but several methods show good empirical performances.

\vspace{2mm}

\noindent $\bullet$ {\bf Deep 2BSDE scheme  \cite{BEJ19}.} 

\vspace{1mm}

\noindent
This scheme  relies on the 2BSDE representation of \cite{CSTV07}
\begin{equation} \label{2BSDE}
\left\{
\begin{array}{ccl}
Y_t & = & g(\Xc_T) - \int_t^T F(s,\Xc_s,Y_s,Z_s,\Gamma_s) \di s - \int_t^T  Z_s\trans \sigma \di W_s, \\
Z_t & = & D_x g(\Xc_T) - \int_t^T A_s \di s - \int_t^T \Gamma_s \sigma  \di W_s, \quad 0 \leq t \leq T, 
\end{array}
\right. 
\end{equation}
with $A_t$ $=$ $\Lc D_x u(t,\Xc_t)$. 
The idea is to adapt the Deep BSDE algorithm to the fully non-linear case. Again, we treat the backward system \eqref{2BSDE} as a forward equation by appro\-ximating the initial conditions $Y_0, Z_0$ and the  $A,\Gamma$ components of the 2BSDE at each time by networks functions taking the forward process
$\Xc$ as input, and aiming to match the terminal condition.

\vspace{2mm}

\noindent $\bullet$ {\bf Second order DBDP (2DBDP) \cite{phawar19}.} 

\vspace{1mm}

\noindent
The basic idea is to adapt the DBDP scheme by  approximating the solution $u$ and its gradient $D_x u$ by network functions $\Uc$ and $\Zc$, and then Hessian  $D_x^2 u$ by the automatic differentiation $D_x\Zc$ of the network function $\Zc$ (or double automatic differentiation $D_x^2\Uc$ of the network function 
$\Uc$), via a learning approach relying on the time discretization of the BSDE \eqref{BSDE2}. 
It turns out that such method approximates poorly $\Gamma$ inducing instability  of the scheme: indeed, while the unique pair solution $(Y,Z)$ to classical BSDEs \eqref{BSDEfeyn} completely characterizes the solution to the related semilinear  PDE and its gradient, the relation  \eqref{BSDE2}  
does not allow to characterize directly the triple $(Y,Z,\Gamma)$. This approach was proposed and tested in \cite{phawar19} where the automatic differentiation is performed on the previous value of $\Zc$ with a truncation $\Tc$ which allows to reduce instabilities.

\vspace{2mm}

\noindent $\bullet$ {\bf Second Order Multistep schemes.} 
\vspace{1mm}

\noindent To overcome the instability in the approximation of the $\Gamma$-component in the Second order DBDP scheme, we propose a finer  approach based on a suitable probabilistic representation of the $\Gamma$-component for learning accurately the Hessian function $D_x^2 u$ by using also  
Malliavin weights. We start from the training simulations of the forward process $(X_i)_i$ on the uniform grid $\pi$ $=$ $\{t_i = i |\pi|,i=0,\ldots,N\}$, $|\pi|$ $=$ $T/N$, and  notice that $X_i$ $=$ $\Xc_{t_i}$, $i$ $=$ $0,\ldots,N$ as $\mu$ and $\sigma$ are constants. 
The approximation of the value function $u$ and its gradient $D_x u$ is learnt simultaneously on the grid $\pi$ but requires in addition 
a preliminary approximation of the Hessian  $D_x^2 u$  in the fully non-linear case. 
This will be performed by regression-based machine learning  scheme on a subgrid $\hat\pi$ $\subset$ $\pi$, which allows to reduce the computational time of the algorithm.

We propose three versions of  second order MDBDP based on different representations of the Hessian function.  
For the second and the third one, we need to introduce a subgrid $\hat\pi$ $=$ $\{t_{\hat\kappa \ell}, \ell =0,\ldots,\hat N\}$ $\subset$ $\pi$, of modulus $|\hat\pi|$ $=$ $\hat\kappa|\pi|$, for some $\hat\kappa$ $\in$ $\N^*$, with $N$ $=$  $\hat\kappa\hat N$.
\begin{itemize}
\item[-] \textit{Version 1:}  Extending the methodology introduced in \cite{phawar19}, the current $\Gamma$-component at step $i$ can be estimated by automatic differentiation of the $Z$-component at the previous step while 
the other $\Gamma$-components are estimated by automatic differentiation  of their associated $Z$-components:
\begin{equation}
\Gamma_i \; \simeq \;  D_x Z_{i+1}, \quad \Gamma_j \; \simeq  \; D_x Z_j, \; j >i. 
\end{equation}
\item[-] \textit{Version 2:}
The time discretization of \eqref{BSDE2} on the time grid $\hat\pi$, where $(Y_\ell^{\hat\pi},Z_\ell^{\hat\pi},\Gamma_\ell^{\hat\pi})$ denotes an approximation of  the triple 
$$
\big(u(t_{\hat\kappa\ell},X_{\hat\kappa\ell}),D_x u(t_{\hat\kappa\ell},X_{\hat\kappa\ell}),D_x^2 u(t_{\hat\kappa\ell},X_{\hat\kappa\ell})\big), \quad \ell = 0,\ldots,\hat N,
$$
leads to the standard  representation formula for the $Z$ component: 
\begin{align}
Z_\ell^{\hat\pi} & = \; \E_{\hat\kappa\ell}  \Big[ Y_{\ell +1}^{\hat\pi}  \hat H_\ell^1 \Big], \quad  \ell = 0,\ldots,\hat N-1,
\end{align}
(recall that $\E_{\hat\kappa \ell}$ denotes the conditional expectation w.r.t. $\mathcal{F}_{t_{\hat\kappa\ell}}$), with the Malliavin weight of order one: 
\begin{align}
\hat H_\ell^1 \; = \; (\sigma\trans)^{-1} \frac{\hat\Delta W_\ell}{|\hat\pi|}, 
\quad  \hat\Delta W_\ell \; := \;  W_{t_{\hat\kappa(\ell+1)}} - W_{t_{\hat\kappa \ell}}.   
\end{align} 
 By direct differentiation, we then obtain an approximation of the $\Gamma$ component as  
\begin{align}
\Gamma_\ell^{\hat\pi} & \simeq  \;  \E_{\hat\kappa\ell}  \Big[ D_x u(t_{\hat\kappa(\ell+1)},X_{\hat\kappa(\ell+1)}) \hat H_\ell^1 \Big].  
\end{align}
Moreover, by introducing the antithetic variable 
\begin{align} \label{antiX1} 
\hat X_{\hat\kappa(\ell+1)}  & = \; X_{\hat\kappa\ell} - \sigma \hat\Delta W_\ell,  
\end{align} 
we then propose the following  regression estimator of $D_x^2u$ on the grid $\hat\pi$ for $\ell = 0,\ldots,\hat N-1$ with
\begin{equation} \label{estim2Gamma}
\left\{
\begin{array}{ccl}
\hat\Gamma^{(1)}(t_{\hat\kappa\hat N},X_{\hat\kappa\hat N})  & = & D^2 g(X_{\hat\kappa\hat N})  \\
\hat\Gamma^{(1)}(t_{\hat\kappa\ell},X_{\hat\kappa\ell})  & = &  
\E_{\hat\kappa\ell}  \Big[ \frac{D_x u(t_{\hat\kappa(\ell+1)},X_{\hat\kappa(\ell+1)}) - D_x u(t_{\hat\kappa(\ell+1)},\hat X_{\hat\kappa(\ell+1)})  }{2}  \hat H_\ell^1 \Big]. 
\end{array}
\right.
\end{equation}
\item[-] \textit{Version 3:} 
 Alternatively, the time discretization of  \eqref{BSDE2} on $\hat\pi$ yields the iterated conditional expectation relation: 
\begin{align}
Y_\ell^{\hat\pi} &= \; \E_{\hat\kappa \ell} \Big[  g(X_{\hat\kappa\hat N}) 
-  |\hat\pi|   \sum_{m=\ell}^{\hat N-1} F(t_{\hat\kappa m},X_{\hat\kappa m},Y_m^{\hat\pi},Z_m^{\hat\pi},\Gamma_m^{\hat\pi})\Big], \quad \ell=0,\ldots,\hat N,  
\end{align}
By (double) integration by parts, and using Malliavin weights on  the Gaussian vector $X$, we obtain a multistep  approximation  of the $\Gamma$-component: 
\begin{align}
\Gamma_\ell^{\hat\pi} & \simeq \;  
\E_{\hat\kappa \ell} \Big[  g(X_{\hat\kappa\hat N}) \hat H_{\ell,\hat N}^2 
- |\hat\pi|   \sum_{m=\ell+1}^{\hat N-1} F(t_{\hat\kappa m},X_{\hat\kappa m},Y_m^{\hat\pi},Z_m^{\hat\pi},\Gamma_m^{\hat\pi})\hat H_{\ell,m}^2  \Big], 
\end{align}
for $\ell = 0,\ldots,\hat N$, where   
 \begin{align}
 \hat H_{\ell,m}^2  & = \;  (\sigma\trans)^{-1} 
 \frac{ \hat\Delta W_\ell^m (\hat\Delta W_\ell^m)\trans - (m-\ell) |\hat\pi| I_d}{(m-\ell)^2|\hat\pi|^2} \sigma^{-1},\  
 \hat\Delta W_\ell^m \; := \;  W_{t_{\hat\kappa m}} - W_{t_{\hat\kappa \ell}}.   
 \end{align}
By introducing again the antithetic variables 
\begin{align} \label{antiX} 
\hat X_{\hat\kappa m}  & = \; X_{\hat\kappa\ell} - \sigma \hat\Delta W_\ell^m, \quad m=\ell+1,\ldots,\hat N, 
\end{align}  
we then propose another regression estimator of $D_x^2u$ on the grid $\hat\pi$  with 
\begin{align}\label{estim1Gamma}
& \hat\Gamma^{(2)}(t_{\hat\kappa\ell},X_{\hat\kappa\ell})  \\ & = \; \E_{\hat\kappa \ell} \Big[  
\frac{g(X_{\hat\kappa\hat N}) + g(\hat X_{\hat\kappa\hat N})}{2} \hat H_{\ell,\hat N}^2  \\
&  \; - \;  \frac{|\hat\pi|}{2}   \sum_{m=\ell+1}^{\hat N-1}  \Big( F\big(t_{\hat\kappa m},X_{\hat\kappa m},u(t_{\hat\kappa m},X_{\hat\kappa m}),
D_x u(t_{\hat\kappa m},X_{\hat\kappa m}),\hat\Gamma^{(2)}(t_{\hat\kappa m},X_{\hat\kappa m})\big) \\
&  + \;  F\big(t_{\hat\kappa m},\hat X_{\hat\kappa m},u(t_{\hat\kappa m},\hat X_{\hat\kappa m}),D_x u(t_{\hat\kappa m},\hat X_{\hat\kappa m}),
\hat\Gamma^{(2)}(t_{\hat\kappa m},\hat X_{\hat\kappa m})\big) \\
&   \; - \;  2  F\big(t_{\hat\kappa \ell},X_{\hat\kappa \ell},u(t_{\hat\kappa \ell},X_{\hat\kappa \ell}),
D_x u(t_{\hat\kappa \ell},X_{\hat\kappa \ell}),\hat\Gamma^{(2)}(t_{\hat\kappa \ell},X_{\hat\kappa \ell})\big)  \Big) \hat H_{\ell,m}^2  \Big],
\end{align}
for $\ell$ $=$ $0,\ldots,N-1$, and $\hat\Gamma^{(2)}(t_{\hat\kappa\hat N},X_{\hat\kappa\hat N})$ $=$ $D^2 g(X_{\hat\kappa\hat N})$.  
The correction term $-2F$ evaluated at time $t_{\hat\kappa\ell}$ in 
$\hat\Gamma^{(2)}(t_{\hat\kappa\ell},X_{\hat\kappa\ell})$ does not add bias since 
\begin{align}
\E_{\hat\kappa \ell} \Big[  F\big(t_{\hat\kappa \ell},X_{\hat\kappa \ell},u(t_{\hat\kappa \ell},X_{\hat\kappa \ell}),
D_x u(t_{\hat\kappa \ell},X_{\hat\kappa \ell}),\hat\Gamma^{(2)}(t_{\hat\kappa \ell},X_{\hat\kappa \ell})\big) \hat H_{\ell,m}^2  \Big] &= \; 0,  
\end{align}
for all $m=\ell+1,\ldots,\hat N-1$, and by Taylor expansion of $F$ at second order, we see that it allows together with the antithetic variable to control the variance when the time step goes to zero.  
\end{itemize}

\vspace{1mm}

\begin{Remark}
{\rm In the case where the function $g$ has some regularity property, one can avoid the integration by parts at the terminal data component in the above expression of $\hat\Gamma^{(2)}$. For example, 
when $g$ is $C^1$, $\frac{g(X_{\hat\kappa\hat N}) + g(\hat X_{\hat\kappa\hat N})}{2} \hat H_{\ell,\hat N}^2$ is alternatively replaced in $\hat\Gamma^{(2)}$ expression by $(Dg(X_{\hat\kappa\hat N}) - Dg(\hat X_{\hat\kappa\hat N}))\hat H_{\ell,\hat N}^1$, while when it is $C^2$ it is replaced by $D^2g(X_{\hat\kappa\hat N})$.
}
\ep
\end{Remark}

\begin{Remark}
{\rm We point out that in our machine learning setting for the versions 2 and 3 of the scheme, we only solve two optimization problems by time step instead of three as in \cite{FTW11}. One optimization is dedicated to the computation of the $\Gamma$ component but the $\Uc$ and $\Zc$ components are simultaneously learned by the algorithm.
}
\ep
\end{Remark}


 We can now describe the three versions of second order MDBDP schemes for the numerical resolution of the fully non-linear PDE  
 \eqref{eq: fully nonlinear PDE}.  We emphasize that these schemes do not require {\it a priori}  that the solution to the PDE is smooth.


\begin{small}

 \begin{algorithm2e}[H]
\DontPrintSemicolon 
\SetAlgoLined 

\For{$i= N -1,\ldots,0$}
    {If $i$ $=$ $N - 1$, update $\widehat{\Gamma}_{i}$ $=$ $D^2 g$, otherwise  $\widehat{\Gamma}_{i}$ $=$ $D_x\widehat{\Zc}_{i+1}$, 
    $\widehat{\Gamma}_{j}$ $=$ $D_x\widehat{\Zc}_{j}$,  $j \in \llbracket i+1,N-1\rrbracket$,  \tcc*{ Update Hessian}
     {Minimize over network functions $\Uc$ $:$ $\R^d$ $\rightarrow$ $\R$, and $\Zc$ $:$ $\R^d$ $\rightarrow$ $\R^d$ the loss function at time 
      $t_i$: 
     \begin{align}
      & J_i^{MB}(\Uc,\Zc)\\ &= \;  \E \Big| g(X_N) - |\pi| \sum_{j=i+1}^{N-1} F(t_j,X_j,\widehat{\Uc}_j(X_j),\widehat{\Zc}_j(X_j),\hat\Gamma_j(X_{j}))  \nonumber
     \\ & \quad \quad \quad - \sum_{j=i+1}^{N-1} \widehat{\Zc}_j(X_j) \trans\sigma  \Delta W_j - \;  \Uc(X_i) \nonumber \\
& \quad \quad \quad  -  |\pi| F(t_i,X_i,\Uc(X_i),\Zc(X_i),\widehat\Gamma_i(X_{i+1}))  -  \Zc(X_i)\cdot \sigma \Delta W_i \Big|^2. 
      \end{align}
     Update    $(\widehat{\Uc}_i,\widehat{\Zc}_i)$ as the solution to this minimization problem  \tcc*{Update the function and its derivative}}
}
\caption{ Second order Explicit Multistep  DBDP (2EMDBDP) \label{algoGlobSchemeExplicit}}
\end{algorithm2e}
 
 
\begin{algorithm2e}
\DontPrintSemicolon 
\SetAlgoLined 
\For{ $\ell$ $=$ $\hat N ,\ldots,0$}
{If $\ell$ $=$ $\hat N$, update $\widehat{\Gamma}_{\ell}$ $=$ $D^2 g$, otherwise
minimize over network functions $\Gamma$ $:$ $\R^d$ $\rightarrow$ $\S^d$ the loss function
\begin{align}
\Jc_\ell^{1,M}(\Gamma) & = \;  \E \Big|  \Gamma(X_{\hat\kappa\ell}) -   \frac{\widehat{\Zc}_{\hat\kappa(\ell+1)}(X_{\hat\kappa(\ell+1)}) - \widehat{\Zc}_{\hat\kappa(\ell+1)}(\hat X_{\hat\kappa(\ell+1)})  }{2}  
\hat H_\ell^1 \Big|^2. 
\end{align}
Update $\widehat\Gamma_\ell$ the solution to this minimization problem  \tcc*{ Update Hessian}
    \For{$k= \hat\kappa -1,\ldots,0$}
     {Minimize over network functions $\Uc$ $:$ $\R^d$ $\rightarrow$ $\R$, and $\Zc$ $:$ $\R^d$ $\rightarrow$ $\R^d$ the loss function at time 
      $t_i$, $i$ $=$ $(\ell-1) \hat\kappa  +k$: 
     \begin{align}
      & J_i^{MB}(\Uc,\Zc)\\ &= \;  \E \Big| g(X_N) - |\pi| \sum_{j=i+1}^{N-1} F(t_j,X_j,\widehat{\Uc}_j(X_j),\widehat{\Zc}_j(X_j),\widehat\Gamma_\ell(X_j)) \nonumber
    \\ & \quad \quad \quad - \sum_{j=i+1}^{N-1} \widehat{\Zc}_j(X_j) \trans\sigma  \Delta W_j - \;  \Uc(X_i) \nonumber \\
& \quad \quad \quad  -  |\pi| F(t_i,X_i,\Uc(X_i),\Zc(X_i),\widehat\Gamma_\ell(X_i))  -  \Zc(X_i)\cdot \sigma \Delta W_i \Big|^2. 
      \end{align}
     Update    $(\widehat{\Uc}_i,\widehat{\Zc}_i)$ as the solution to this minimization problem  \tcc*{Update the function and its derivative}
    }
}
\caption{ Second order Multistep  DBDP (2MDBDP) \label{algoGlobScheme}}
\end{algorithm2e}


\begin{algorithm2e}
\DontPrintSemicolon 
\SetAlgoLined 
\For{ $\ell$ $=$ $\hat N ,\ldots,0$}
{If $\ell$ $=$ $\hat N$, update $\widehat{\Gamma}_{\ell}$ $=$ $D^2 g$, otherwise
minimize over network functions $\Gamma$ $:$ $\R^d$ $\rightarrow$ $\S^d$ the loss function
\begin{align}
& \Jc_\ell^{2,M}(\Gamma)\\ & = \;  \E \Big|  \Gamma(X_{\hat\kappa\ell}) - 
\frac{D^2g(X_{\hat\kappa\hat N})+D^2g(\hat X_{\hat\kappa\hat N})}{2}\\
& \quad \quad \; +  \;  \frac{|\hat\pi|}{2}   \sum_{m=\ell+1}^{\hat N-1}  \Big( F\big(t_{\hat\kappa m},X_{\hat\kappa m}, \widehat{\Uc}_{\hat\kappa m}(X_{\hat\kappa m}),
\widehat{\Zc}_{\hat\kappa m}(X_{\hat\kappa m}),\widehat\Gamma_{m}^{}(X_{\hat\kappa m})\big) \\
& \quad \quad \quad \quad  \quad + \;  F\big(t_{\hat\kappa m},\hat X_{\hat\kappa m},\widehat{\Uc}_{\hat\kappa m}(\hat X_{\hat\kappa m}),\widehat{\Zc}_{\hat\kappa m}(\hat X_{\hat\kappa m}),
\widehat{\Gamma}^{}_{m}(\hat X_{\hat\kappa m})\big)  \\
&  \quad \quad \quad \quad  \quad -  \; 2  F\big(t_{\hat\kappa \ell},\hat X_{\hat\kappa \ell},\widehat{\Uc}_{\hat\kappa \ell}(\hat X_{\hat\kappa \ell}),\widehat{\Zc}_{\hat\kappa \ell}(\hat X_{\hat\kappa \ell}),
\widehat{\Gamma}^{}_{\ell}(\hat X_{\hat\kappa \ell})\big)   \Big) \hat H_{\ell,m}^2  \Big|^2. 
\end{align}
Update $\widehat{\Gamma}_\ell$ the solution to this minimization problem  \tcc*{ Update Hessian}
    \For{$k= \hat\kappa -1,\ldots,0$}
     {Minimize over network functions $\Uc$ $:$ $\R^d$ $\rightarrow$ $\R$, and $\Zc$ $:$ $\R^d$ $\rightarrow$ $\R^d$ the loss function at time 
      $t_i$, $i$ $=$ $(\ell-1) \hat\kappa  +k$: 
     \begin{align}
      & J_i^{MB}(\Uc,\Zc)\\ &= \;  \E \Big| g(X_N) - |\pi| \sum_{j=i+1}^{N-1} F(t_j,X_j,\widehat{\Uc}_j(X_j),\widehat{\Zc}_j(X_j),\widehat{\Gamma}_\ell(X_j)) \nonumber 
     \\ & \quad \quad \quad - \sum_{j=i+1}^{N-1} \widehat{\Zc}_j(X_j) \trans\sigma  \Delta W_j  - \;  \Uc(X_i) \nonumber \\
& \quad \quad \quad  -  |\pi| F(t_i,X_i,\Uc(X_i),\Zc(X_i),\widehat{\Gamma}_\ell(X_i))  -  \Zc(X_i)\cdot \sigma \Delta W_i \Big|^2.  
      \end{align}
     Update    $(\widehat{\Uc}_i,\widehat{\Zc}_i)$ as the solution to this minimization problem  \tcc*{Update the function and its derivative}
    }
}
\caption{ Second order Multistep Malliavin DBDP (2M$^2$DBDP) \label{algoGlobScheme2}}
\end{algorithm2e}

\end{small}


The  proposed algorithms \ref{algoGlobSchemeExplicit}, \ref{algoGlobScheme}, \ref{algoGlobScheme2}  are in  backward iteration, and involve one  optimization at each step. 
Moreover, as the computation of   $\Gamma$ requires a further derivation for Algorithms \ref{algoGlobScheme} and \ref{algoGlobScheme2}, 
we may expect that  the additional propagation error varies according to  $\frac{|\pi|}{|\hat\pi|}$ $=$ $\frac{1}{\hat\kappa}$, and thus the convergence of the scheme when $\hat\kappa$ is large. 
In the numerical implementation, the expectation in the loss functions 
are replaced by empirical average and the minimization over network functions is performed by stochastic gradient descent. 

\section{Numerical applications}\label{sec: numerics}

We test our different algorithms on various examples and by varying  the state space dimension. 
If not stated otherwise, we choose the maturity $T=1$. 
In each example we use an architecture composed of 2 hidden layers with $d+10$ neurons. We apply Adam gradient descent \cite{KB14} with a decreasing learning rate, using the Tensorflow library \cite{tensorflow}. 
Each numerical experiment is conducted using a node composed of 2 Intel® Xeon® Gold 5122 Processors, 192 Go of RAM, and 2 GPU nVidia® Tesla® V100 16Go.  We use a batch size of 1000.


\subsection{Numerical tests on credit valuation adjustment pricing}

We consider an example of model from \cite{labordere} for the pricing of CVA in a $d$-dimensional Black-Scholes model
\begin{equation}
\di X_t = \sigma X_t\ \di W_t,\ X_0 = 1_d
\end{equation} with $\sigma > 0$, given by the nonlinear PDE
\begin{align}
\begin{cases}
     \partial_t u + \frac{\sigma^2 }{2} \mathrm{Tr}(x^\top D^2_x u\ x) + \beta (u_+ - u) = 0\ & \mathrm{on}\ [0,T]\times \R^d\\
     u(T,x) =  |\sum_{i=1}^d x_i - d| - 0.1 & \mathrm{on}\ \R^d
    \end{cases}
\end{align} with a straddle type payoff.
We compare our results  with the DBDP scheme \cite{HPW19} with the ones from the  Deep BSDE solver \cite{HJE17}.  The results in Table \ref{tab: CVA}  are  averaged over 10 runs and the standard deviation is written in parentheses. We use ReLu activation functions.

\begin{table}[H]
    \centering
    \begin{tabular}{|c|c|c|}
    \hline
   Dimension $d$  & DBDP \cite{HPW19}  & DBSDE \cite{HJE17}  \\
    \hline
    1 & 0.05950 (0.000257) & 0.05949 (0.000264)
\\
    \hline      
    3 & 0.17797 (0.000421) & 0.17807 (0.000288) \\
    \hline
    5 & 0.25956 (0.000467) & 0.25984 (0.000331)\\
    \hline
    10 & 0.40930 (0.000623) & 0.40886 (0.000196)
\\
    \hline
    15 & 0.52353 (0.000591) & 0.52389 (0.000551)\\
    \hline
    30 & 0.78239 (0.000832) & 0.78231 (0.001266)\\
    \hline
    \end{tabular}
    \caption{CVA value with $X_0 = 1, T = 1,\beta = 0.03, \sigma = 0.2$ and $50$ time steps. }
    \label{tab: CVA}
\end{table}
We observe in Table \ref{tab: CVA} that both algorithms give very close results and are able to solve the nonlinear pricing problem in high dimension $d$. The variance of the results is quite small and similar from one to another but increases with the dimension. The same conclusions arise when solving the PDE for the larger maturity $T=2$.

\subsection{Portfolio allocation in stochastic volatility models}

We consider several examples from \cite{phawar19} that we solve with Algorithms \ref{algoGlobSchemeExplicit} (2EMDBDP), 
\ref{algoGlobScheme} (2MDBDP), and \ref{algoGlobScheme2} (2M$^2$DBDP) designed in this paper.  Notice that some comparison tests with the 2DBSDE scheme 
\cite{BEJ19} have been already done in \cite{phawar19}. For a resolution with $N = 120, {\hat N} = 30$, the execution of our multitep algorithms takes between 10000 s. and 30000 s. (depending on the dimension) with a number of gradient descent iterations fixed at 4000 at each time step except 80000 at the first one. We use tanh as activation function. 

We consider a portfolio selection problem formulated as follows. There are $n$ risky assets of uncorre\-lated price process $P$ $=$ $(P^1,\ldots,P^n)$ with dynamics governed by 
\begin{align}
\di P_t^i  & = \; P_t^i  \sigma(V_t^i)   \big[ \lambda_i (V_t^i)  \di t +   \di W_t^i\big], \quad i=1,\ldots,n,
\end{align} 
where $W$ $=$ $(W^1,\ldots,W^n)$  is a $n$-dimensional Brownian motion,  $\lambda$ $=$ $(\lambda^1,\ldots,\lambda^n)$  is the market price of risk of the assets, $\sigma$ is a positive function (e.g. $\sigma(v)$ $=$ $e^v$ corresponding to the Scott model),  
and  $V$ $=$ $(V^1,\ldots,V^n)$  is the  volatility factor modeled by an Ornstein-Uhlenbeck (O.U.) process
\begin{align} \label{dynY} 
 \di V_t^i  &= \;  \kappa_i[ \theta_i - V_t^i] \di t + \nu_i   \di B^i_t, \quad i=1,\ldots,n, 
\end{align} 
with  $\kappa_i, \theta_i, \nu_i$ $>$ $0$, and 
$B$ $=$ $(B^1,\ldots,B^n)$ a $n$-dimensional Brownian motion, s.t. $d<W^i,B^j>$ $=$ $\delta_{ij} \rho_{ij} dt$, with $\rho_{i}$ $:=$ $\rho_{ii}$ $\in$ $(-1,1)$.  An agent can invest at any time an amount $\alpha_t$ $=$ $(\alpha_t^1,\ldots,\alpha_t^n)$  in the stocks, 
which generates a wealth process $\Xc$ $=$ $\Xc^\alpha$ governed by 
\begin{eqnarray*}
\di \Xc_t  &= &   \sum_{i=1}^n  \alpha_t^i  \sigma(V_t^i) \big[ \lambda_i(V_t^i) \di t +   \di W_t^i \big].  
\end{eqnarray*}
The objective of the agent is to  maximize her expected utility from terminal wealth: 
\begin{align}
\E \big[ U(\Xc_T^\alpha) ] \quad  & \leftarrow  \quad \mbox{maximize over } \alpha 
\end{align}
It is well-known that the solution to this problem can be characterized by the dynamic programming method (see e.g. \cite{pha09}), which leads to the Hamilton-Jacobi-Bellman for the value function on 
$[0,T)\times\R\times\R^n$:  
\begin{equation*}
    \begin{cases}
                \partial_t u  + \Sum_{i=1}^n \big[ \kappa_i(\theta_i - v_i) \partial_{v_i} u + \frac{1}{2} \nu_i^2  \partial_{v_i}^2  u  \big]\\ \; = \; 
       \frac{1}{2} R(v)  \frac{(\partial_\mrx u)^2}{  \partial_{\mrx\mrx}^2 u}  + \sum_{i=1}^n \big[  \rho_i \lambda_i(v_i) \nu_i \frac{\partial_\mrx u \partial_{\mrx v_i}^2 u}{ \partial_{\mrx\mrx}^2 u}  + 
        \frac{1}{2}\rho_i^2 \nu_i^2  \frac{(\partial^2_{\mrx v_i} u)^2}{ \partial^2_{\mrx\mrx} u} \big] 
         &
         \\
        u(T,\mrx,v) \; = \;  U(\mrx),   \quad \quad \mrx \in \R,  \;  v \in \R^n, &
    \end{cases}
\end{equation*} 
with a Sharpe ratio $R(v)$ $:=$ $|\lambda(v)|^2$, for $v$ $=$ $(v_1,\ldots,v_n)$ $\in$ $(0,\infty)^n$.  
The optimal portfolio strategy is then  given in feedback form by 
$\alpha_t^*$ $=$ $\hat a(t,\Xc_t^*,V_t)$, where $\hat a$ $=$ $(\hat a_1,\ldots,\hat a_n)$ is given by 
\begin{align}
& \hat a_i(t,\mrx,v)\\  &= \;   - \frac{1}{\sigma(v_i)} \Big(  \lambda_i(v_i)   \frac{\partial_\mrx u}{ \partial_{\mrx\mrx}^2 u} +   \rho_i \nu_i \frac{\partial_{\mrx v_i}^2 u}{ \partial_{\mrx\mrx}^2 u } \Big),\ (t,\mrx,v=(v_1,\ldots,v_n)) \in [0,T)\times\R\times\R^n, 
\end{align}
for $i$ $=$ $1,\ldots,n$.

\vspace{2mm}

We shall test this example when the utility function $U$ is of exponential form: $U(\mrx)$ $=$ $-\exp(-\eta \mrx)$, with $\eta$ $>$ $0$,  
and under different cases for which  explicit solutions are available. We refer to \cite{phawar19} where these solutions are described.
\begin{itemize}
\item[(1)]  {\it Merton problem.} This corresponds to a degenerate case where the factor $V$, hence the volatility $\sigma$ and  the risk premium $\lambda$  are  constant ($v_i = \theta_i,\ \nu_i = 0$). We train our algorithms with the forward process 
\begin{align}
X_{k+1} & = \; X_k +  |\lambda|  \Delta t_k +   \Delta W_k, \quad k=0,\ldots,N, \; X_0  \;  = \; x_0. 
\end{align} 

\item[(2)]  {\it One risky asset: $n$ $=$ $1$}.  
We train our algorithms  with the forward process 
\begin{align}
\Xc_{k+1} & = \; \Xc_k +  \lambda(\theta)   \Delta t_k +  \Delta W_k,   \quad k=0,\ldots,N-1, \;\;  \Xc_0 \;  = \; \mrx_0 \\
V_{k+1} & = \; V_k  +  \nu \Delta B_k, \;\;\;   \quad k=0,\ldots,N-1, \;\; V_0 \;  = \; \theta.
\end{align}
We test our algorithm with $\lambda(v)$ $=$ $\lambda v$, $\lambda$ $>$ $0$, for which we have an explicit solution.

\item[(3)] {\it No leverage effect, i.e.,  $\rho_i$ $=$ $0$, $i$ $=$ $1,\ldots,n$}.  We train with the forward process
\begin{align}
\Xc_{k+1} & = \; \Xc_k + \Sum_{i=1}^n \lambda_i(\theta_i) \Delta t_k +  \Delta W_k,   \quad k=0,\ldots,N-1, \;\;   \Xc_0 \;  = \; \mrx_0 \\
V_{k+1}^i & = \; V_k^i  +  \nu_i \Delta B_k^i, \;\;\;   \quad k=0,\ldots,N-1, \;\; V_0^i  \;  = \; \theta_i.
\end{align}
  We test our algorithm with $\lambda_i(v)$ $=$ $\lambda_i v_i$, $\lambda_i$ $>$ $0$, $i$ $=$ $1,\ldots,n$, $v$ 
$=$ $(v_1,\ldots,v_n)$,  for which we have an explicit solution.
\end{itemize}

\vspace{2mm}

\begin{figure}[H]
   \begin{minipage}[c]{.49\linewidth}
          \includegraphics[width=0.9\linewidth]{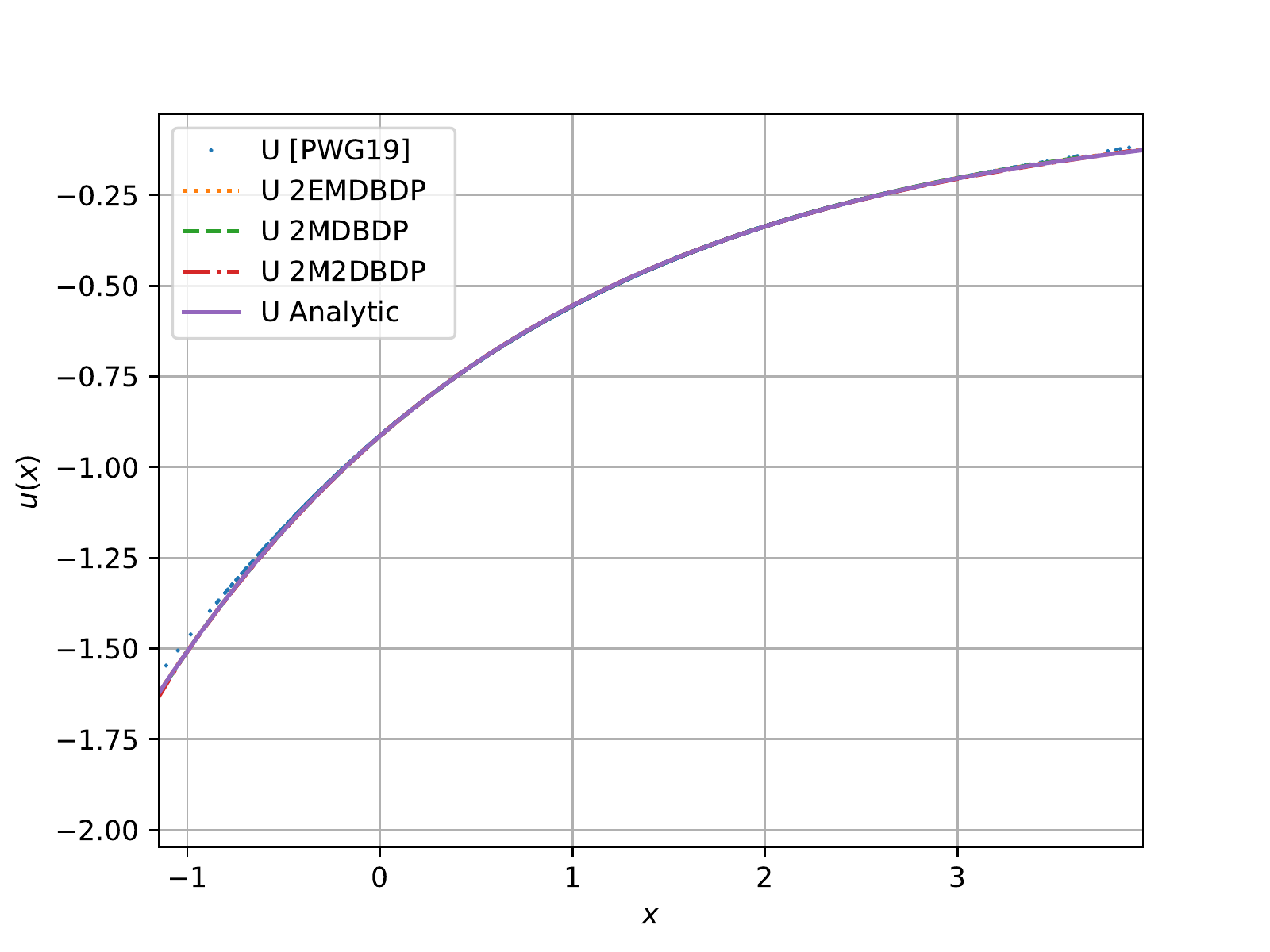}
   \end{minipage} \hfill
   \begin{minipage}[c]{.49\linewidth}
      \includegraphics[width=0.9\linewidth]{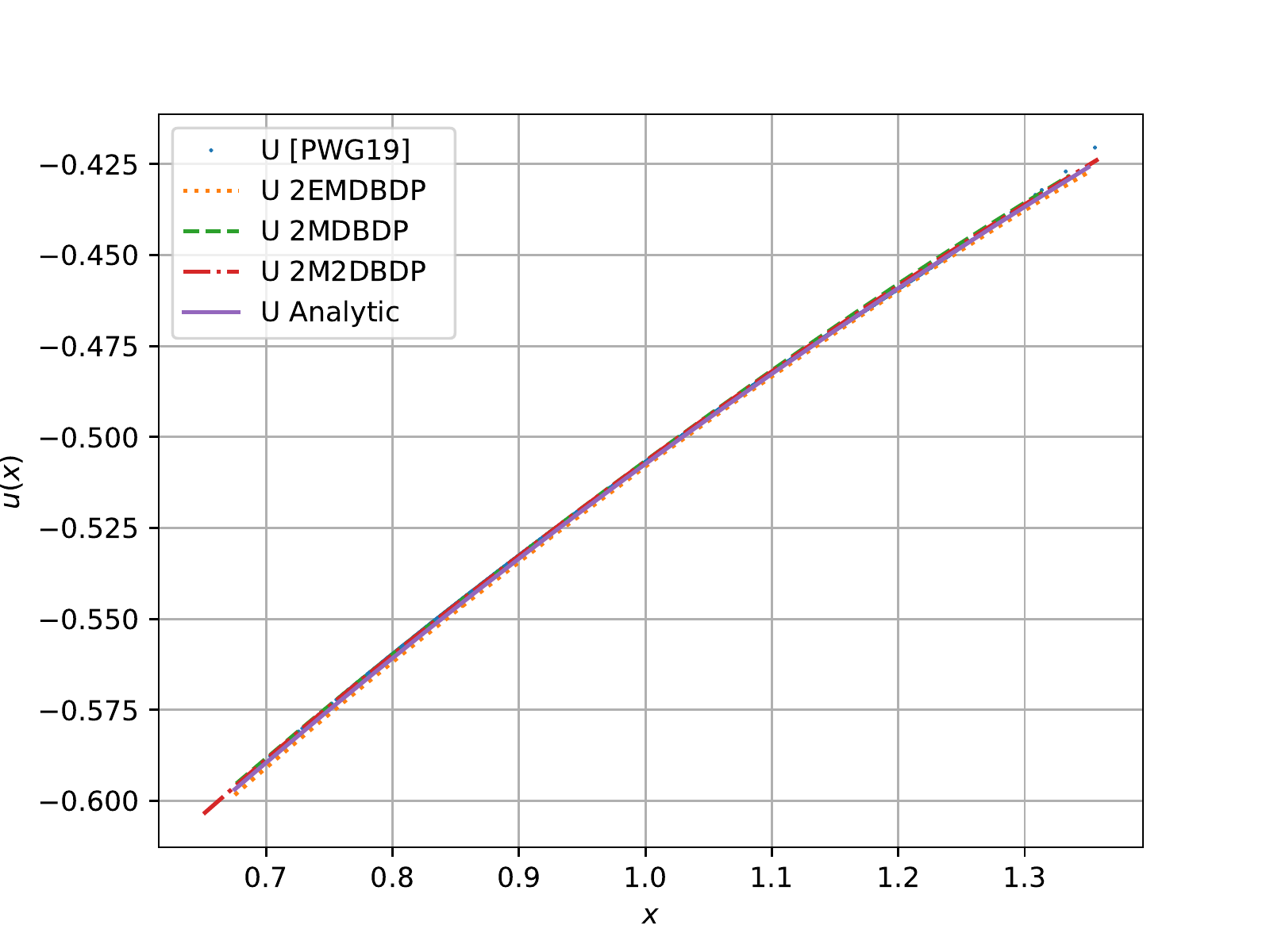}
   \end{minipage}
   \begin{minipage}[c]{.49\linewidth}
          \includegraphics[width=0.9\linewidth]{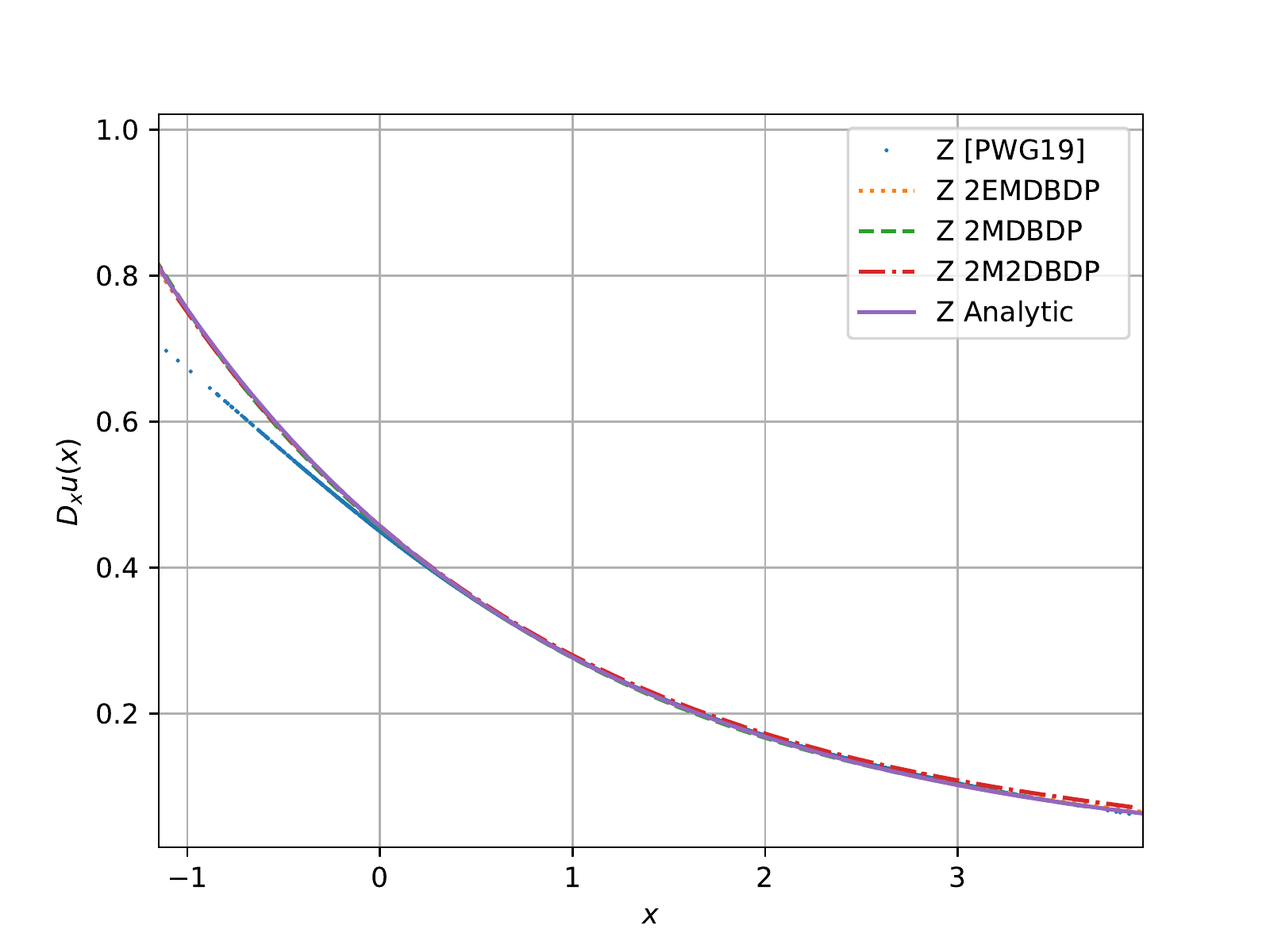}
   \end{minipage} \hfill
   \begin{minipage}[c]{.49\linewidth}
      \includegraphics[width=0.9\linewidth]{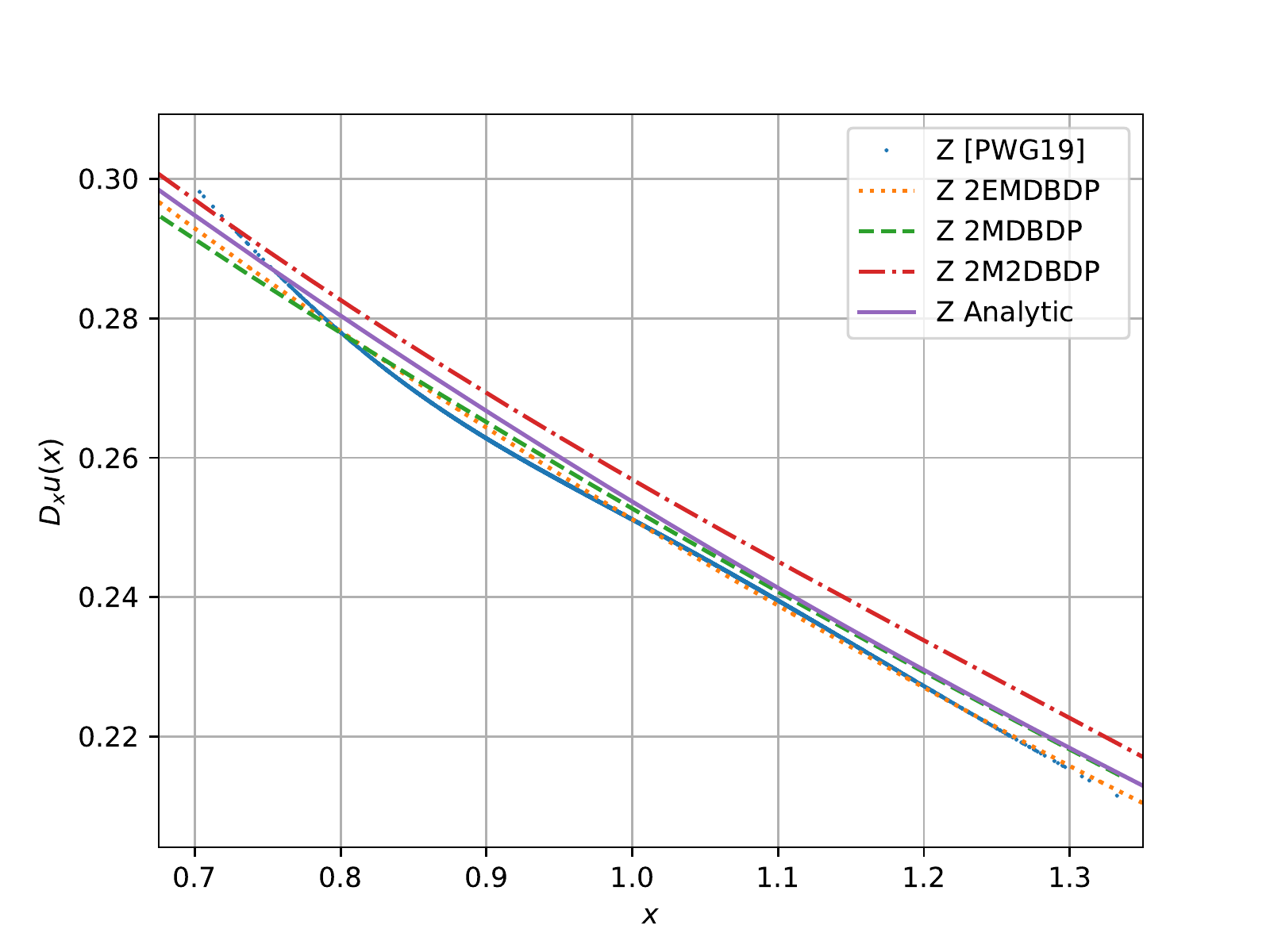}
   \end{minipage}
   \begin{minipage}[c]{.49\linewidth}
          \includegraphics[width=0.9\linewidth]{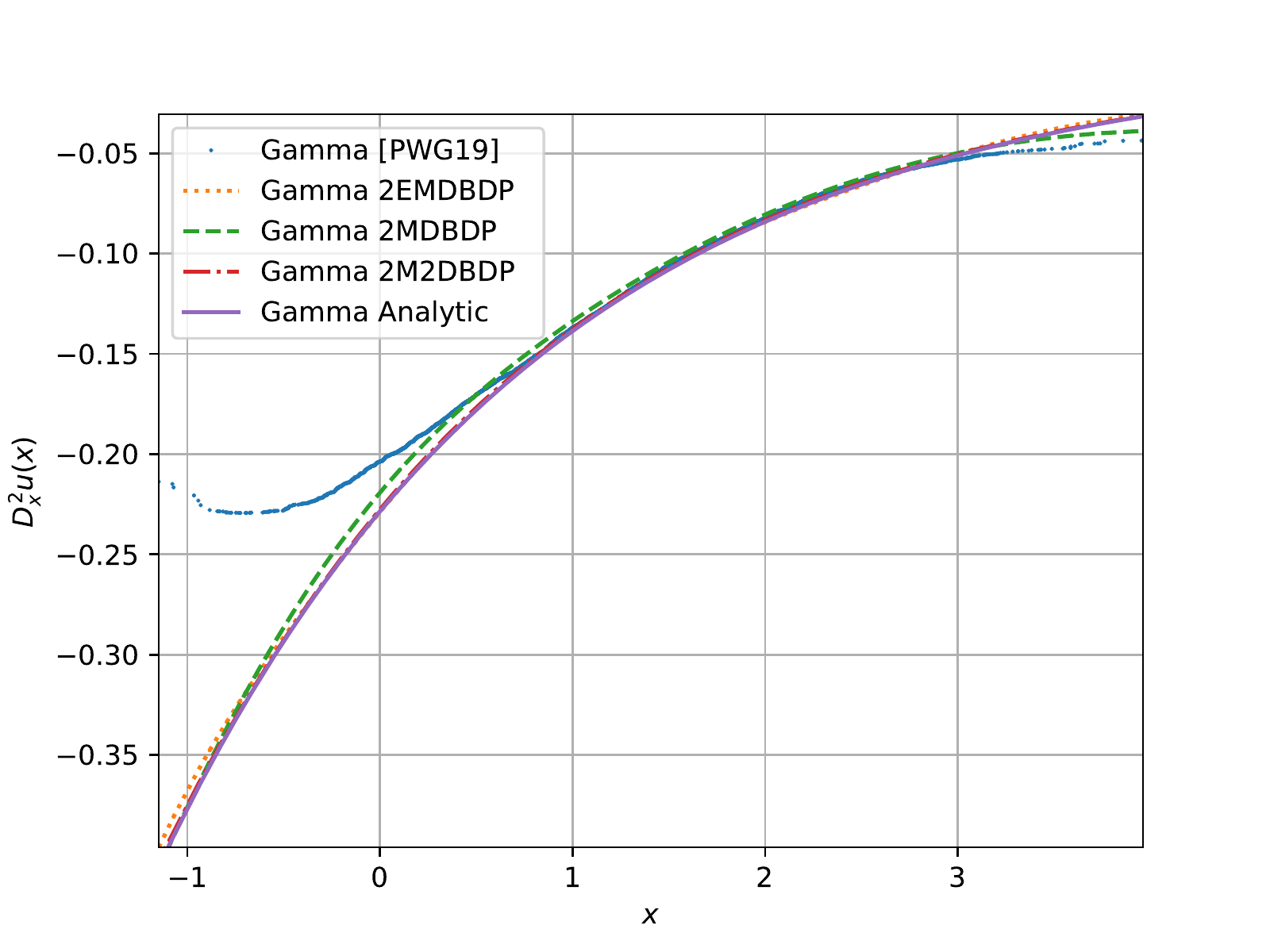}
   \end{minipage} \hfill
   \begin{minipage}[c]{.49\linewidth}
      \includegraphics[width=0.9\linewidth]{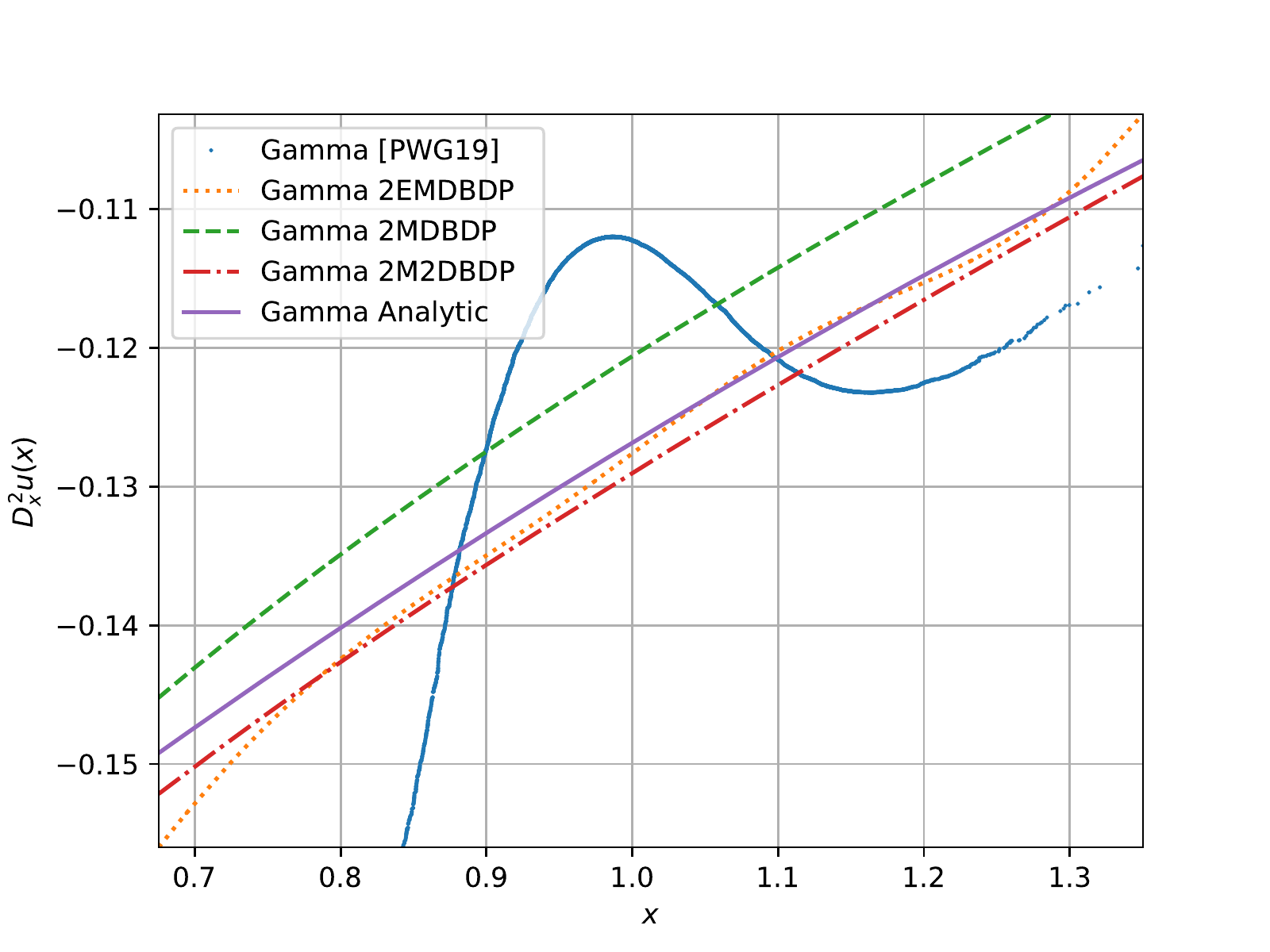}
   \end{minipage}
   \begin{minipage}[c]{.49\linewidth}
          \includegraphics[width=0.9\linewidth]{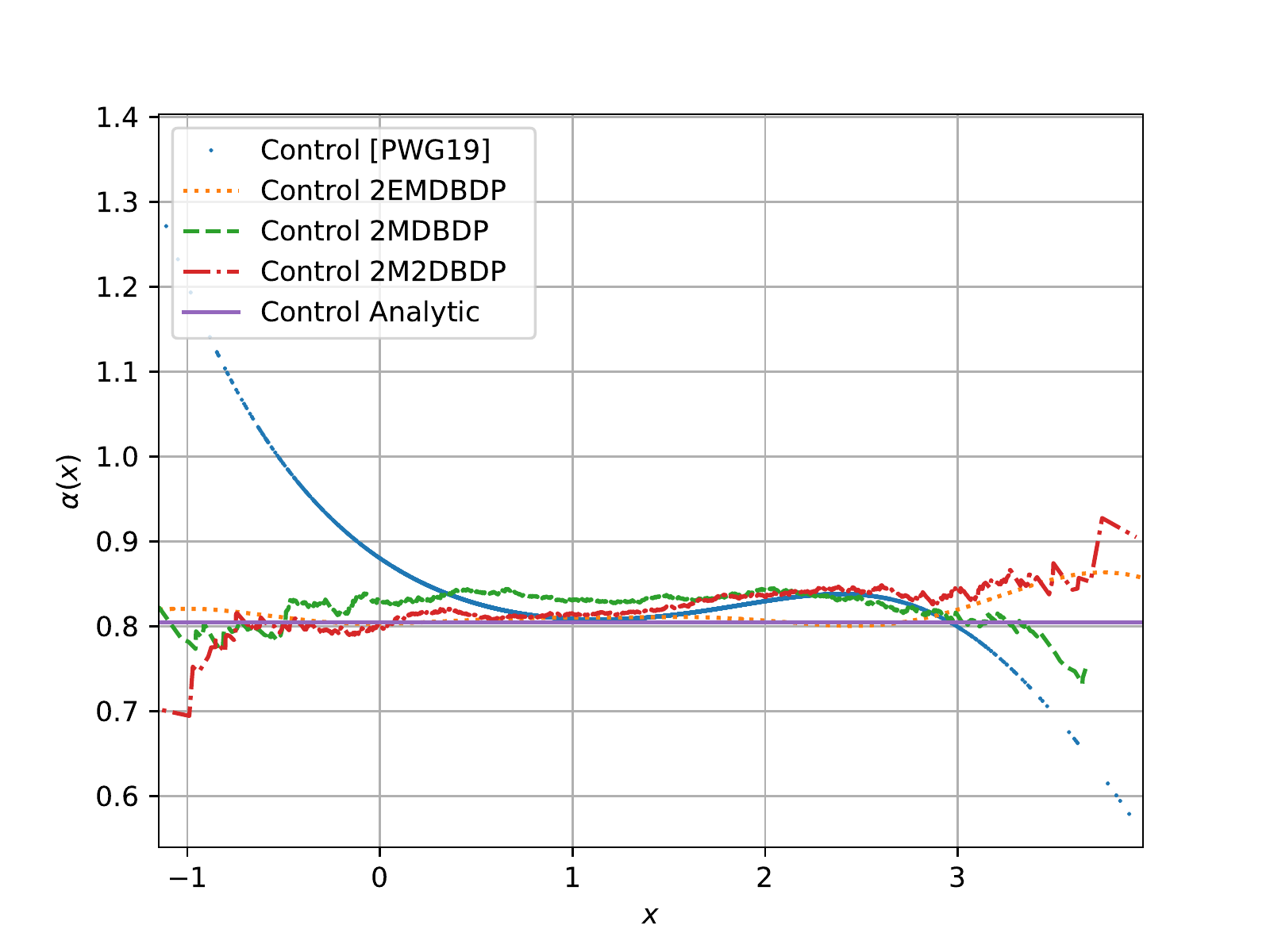}
   \end{minipage} \hfill
   \begin{minipage}[c]{.49\linewidth}
      \includegraphics[width=0.9\linewidth]{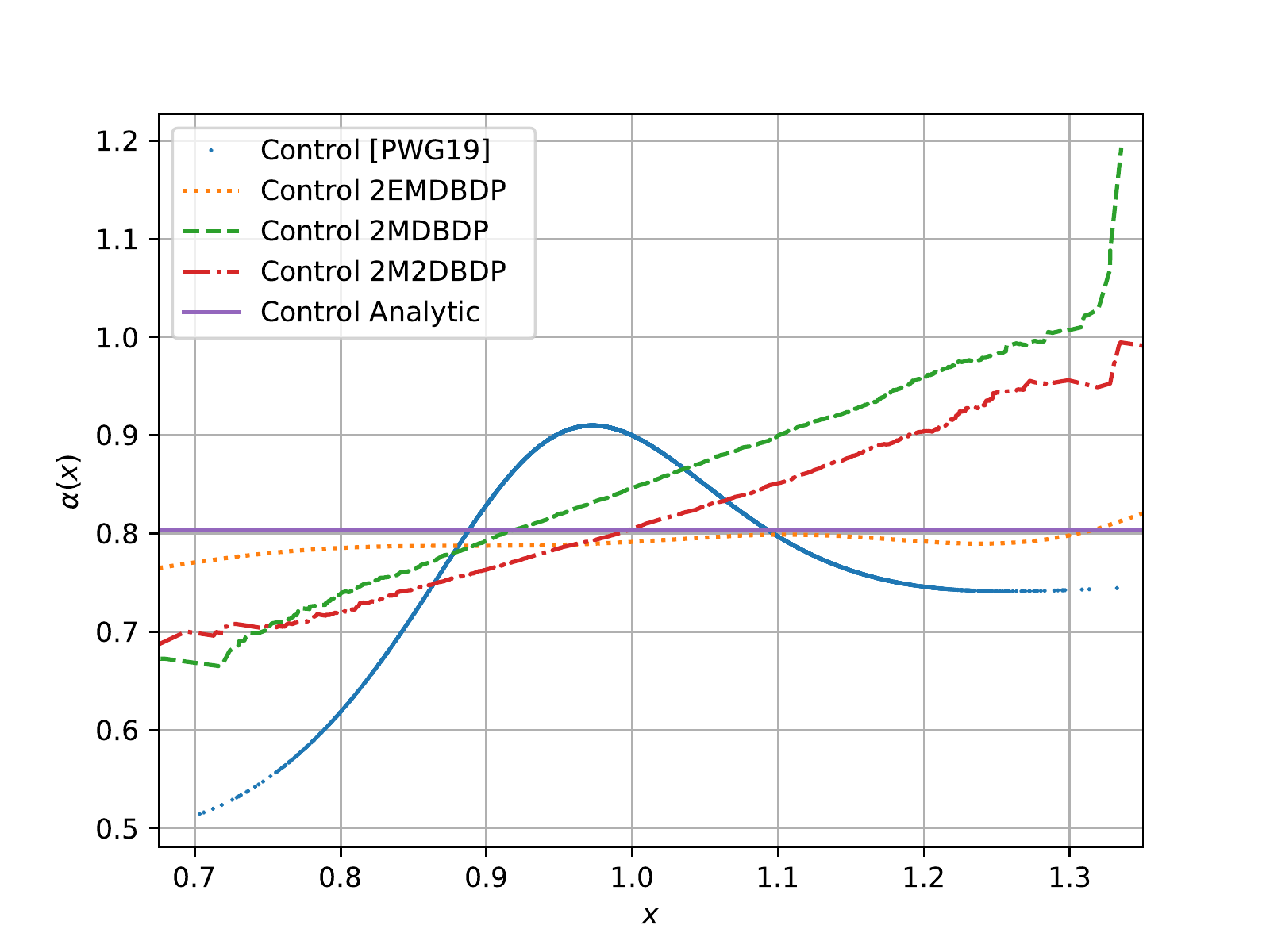}
   \end{minipage}
   \caption{Estimates of $u$, $D_x u$, $D_x^2 u$ and of the optimal control $\alpha$ on the Merton problem  with $N = 120$, $\hat N = 30$. We take $x_0 = 1.$, at the left $t= 0.5042 $, and at the right $t= 0.0084$.} \label{figMerton}
\end{figure}

\textbf{Merton Problem.} 
We take $\eta = 0.5$, $\lambda = 0.6$, $N = 120$, $\hat N = 30$, $T=1$, $x_0 = 1.$
We plot in Figure \ref{figMerton}  the neural networks approximation  
of $u, D_x u, D^2_x u$, and the feedback control $\hat a$ (for one asset) computed  from our different algorithms,  together with their analytic values (in orange). 
As also reported in the estimates of  Table  \ref{fig: table results Merton}, the multistep algorithms improve significantly the results obtained in \cite{phawar19}, where the estimation of the Hessian is not really  accurate (see blue curve in Figure \ref{figMerton}). 


\begin{table}[H]
	\centering
	\begin{tabular}{|c|c|c|c|}
		\hline
		& Average & Standard deviation & Relative error (\%)\\
		\hline
		\cite{phawar19} & 	 -0.50561
 & 0.00029 &  0.20 \\
		\hline
		2EMDBDP & \textbf{-0.50673} & \textbf{0.00019} & \textbf{0.022}\\
		\hline
		2MDBDP & -0.50647 & 0.00033 & 0.030 \\
		\hline
		2M$^2$DBDP & -0.50644 & 0.00022 & 0.035\\
		\hline
	\end{tabular}
	\caption{Estimate of $u(0, 1. )$ in the Merton problem with $N = 120$, $\hat N = 30$. Average and standard deviation observed over 10 independent runs are reported. The theoretical solution is -0.50662. }
	\label{fig: table results Merton}
\end{table}

\textbf{One asset $n$ $=$ $1$ in Scott volatility model.}  
We take $\eta = 0.5$, $\lambda = 1.5$, $\theta = 0.4$, $\nu = 0.4$, $\kappa = 1$, $\rho = -0.7 $, $T=1$, $x_0 = 1.$ For all tests we choose $N = 120$, $\hat N = 30$ and 
$\sigma(v) = e^v$. We  report in Table \ref{fig: table results One asset d = 2} 
the relative error between the neural networks approximation  of $u, D_\mrx u, D^2_\mrx u$ computed from our different algorithms 
and their analytic values. It turns out that the multistep extension of \cite{phawar19}, namely 2EMDBDP scheme,  
yields a very accurate  approximation result, much better than the other algorithms, with also a reduction of the standard deviation.

\begin{table}[H]
	\centering
	\begin{tabular}{|c|c|c|c|}
		\hline
		& Average & Standard deviation & Relative error (\%)\\
		\hline
		\cite{phawar19} & -0.53431 & 0.00070  & 0.34\\
		\hline
		2EMDBDP & \textbf{-0.53613} & \textbf{0.00045}  & \textbf{0.007}\\
		\hline
		2MDBDP & -0.53772
 & 0.00046 & 0.304\\
		\hline
		2M$^2$DBDP & -0.53205  & 0.00050 & 0.755\\
		\hline
	\end{tabular}       
	\caption{Estimate of $u(0, 1, \theta)$ on the One Asset problem with stochastic volatility  ($d=2$) and $N = 120$, $\hat N = 30$. Average and standard deviation observed over 10 independent runs are reported. The exact solution is $-0.53609477$.}
	\label{fig: table results One asset d = 2}
\end{table}

\textbf{No Leverage  in Scott model.} 
In the case with one asset we take $\eta = 0.5$, $\lambda = 1.5$, $\theta = 0.4$, $\nu = 0.2$, $\kappa = 1$, $T=1$, $x_0 = 1.$ For all tests we choose $N = 120$, $\hat N = 30$ and $\sigma(v) = e^v$. We  report in Table \ref{fig: table results No Leverage d = 2} 
the relative error between the neural networks approximation  of $u, D_\mrx u, D^2_\mrx u$ computed  from our different algorithms 
and their analytic values.  All the algorithms yield quite accurate results, but compared to the case with correlation in Table \ref{fig: table results One asset d = 2}, it appears here that the best performance in terms of precision  is achieved by Algorithm 2M$^2$DBDP.

\begin{table}[H]
	\centering
	\begin{tabular}{|c|c|c|c|}
		\hline
		& Average & Standard deviation & Relative error (\%)\\
		\hline
		\cite{phawar19} & -0.49980 & 0.00073 & 0.35\\
		\hline
		2EMDBDP & -0.50400
 & 0.00229 & 0.485\\
		\hline
		2MDBDP &  -0.50149  &  \textbf{0.00024} & 0.015\\
		\hline
		2M$^2$DBDP & \textbf{-0.50157} & 0.00036 & \textbf{0.001}\\
		\hline
	\end{tabular}
	\caption{Estimate of $u(0, 1, \theta)$, with 120 time steps on the No Leverage problem with 1 asset ($d=2$) and $N = 120$, $\hat N = 30$. Average and standard deviation observed over 10 independent runs are reported. The exact solution is $-0.501566 $. }
	\label{fig: table results No Leverage d = 2}
\end{table}  

In the case with four assets ($n$ $=$ $4$, $d=5$), we take $\eta = 0.5$,\\ $\lambda = \begin{pmatrix} 1.5 & 1.1 & 2. & 0.8 \end{pmatrix}$, $\theta = \begin{pmatrix} 0.1 & 0.2 & 0.3 & 0.4 \end{pmatrix}$, 
$\nu = \begin{pmatrix}0.2 & 0.15 & 0.25 & 0.31 \end{pmatrix}$, $\kappa = \begin{pmatrix}1. & 0.8 & 1.1 & 1.3\end{pmatrix}$.
The results are reported in Table \ref{fig: table results No Leverage d = 5}. We observe that the algorithm in \cite{phawar19} provides a not so accurate outcome, while its multistep version (2EMDBDP scheme) divides  by $10$ the relative error and the standard deviation.

\begin{table}[H]
	\centering
	\begin{tabular}{|c|c|c|c|}
		\hline
		& Average & Standard deviation & Relative error (\%)\\
		\hline
		\cite{phawar19} & -0.43768 & 0.00137 & 0.92\\
		\hline
		2EMDBDP & \textbf{-0.4401}
 & \textbf{0.00051} & \textbf{0.239}\\
		\hline
		2MDBDP & -0.43796
 & 0.00098 & 0.861 \\
		\hline
		2M$^2$DBDP & -0.44831
 & 0.00566 & 1.481 \\
		\hline
	\end{tabular}
	\caption{Estimate of $u(0, 1, \theta)$, with 120 time steps on the No Leverage problem  with 4 assets ($d=5$) and $N = 120$, $\hat N = 30$. Average and standard deviation observed over 10 independent runs are reported. The theoretical solution is -0.44176462. }
	\label{fig: table results No Leverage d = 5}
\end{table}  

In the case with nine assets ($n$ $=$ $9$, $d=10$), we take $\eta = 0.5$,\\ $\lambda = \begin{pmatrix} 1.5 & 1.1 & 2. & 0.8 & 0.5 & 1.7 & 0.9 & 1. & 0.9\end{pmatrix}$,\\ $\theta = \begin{pmatrix} 0.1 & 0.2 & 0.3 & 0.4 & 0.25 & 0.15 & 0.18 & 0.08 & 0.91 \end{pmatrix}$,\\ $\nu = \begin{pmatrix}0.2 & 0.15 & 0.25 & 0.31 & 0.4 & 0.35 & 0.22 & 0.4 & 0.15\end{pmatrix}$,\\ $\kappa = \begin{pmatrix}1. & 0.8 & 1.1 & 1.3 & 0.95 & 0.99 & 1.02 & 1.06 & 1.6\end{pmatrix}$.
The results are reported in Table \ref{fig: table results No Leverage d = 10}. 
The approximation is less accurate than in lower dimension, but we observe again that compared to one-step scheme in \cite{phawar19} , 
the multistep versions improve significantly the standard deviation of the result. However the best performance in precision is obtained here by the \cite{phawar19} scheme.

\begin{table}[H]
	\centering
	\begin{tabular}{|c|c|c|c|c|}
		\hline
		 & $\hat N$ & Average & S.d. & Relative error (\%)\\
		\hline
		\cite{phawar19} &  &  \textbf{-0.27920} & 0.05734 & \textbf{1.49}\\
		\hline
		2EMDBDP & & -0.26631
 & \textbf{0.00283} & 3.19\\
		\hline
		2MDBDP & 30 & -0.28979 & 0.00559 & 5.34\\
		\hline
		2MDBDP & 60  & -0.28549
 & 0.00948 & 3.78\\
		\hline
				2MDBDP & 120  & \textbf{-0.28300}
 & 0.01129 & \textbf{2.87}\\
		\hline
		2M$^2$DBDP& 30  & NC & NC & NC \\
		\hline
	\end{tabular}
	\caption{Estimate of $u(0, 1, \theta)$, with 120 time steps on the No Leverage problem with 9 assets ($d=10$) and $N = 120$. Average and standard deviation (S.d.) observed over 10 independent runs are reported. The theoretical solution is -0.27509173.}
	\label{fig: table results No Leverage d = 10}
\end{table}

\section{Extensions and perspectives}\label{sec: extensions}

\noindent $\bullet$ {\bf Solving mean-field control and mean-field games through McKean-Vlasov FBSDEs.} 

\vspace{1mm}

\noindent
These methods solve the optimality conditions for mean-field problems through the stochastic Pontryagin principle from \cite{cardel}. The law of the solution influences the coupled FBSDEs dynamics so they are of McKean-Vlasov type. 
Variations around the Deep BSDE method \cite{HJE17} are used to solve such a system by \cite{carlau19}, \cite{fouquezhang19}. \cite{germicwar19} uses the Merged method from \cite{CWNMW19} and solves several numerical examples in dimension 10 by introducing an efficient law estimation technique. \cite{carlau19} also proposes another method dedicated to mean field control to directly tackle the optimization problem with a neural network as the control in the stochastic dynamics. The $N-$player games, before going to the mean-field limit of an infinite number of players, are solved by \cite{H19}, \cite{HHL19}. 

\vspace{2mm}

\noindent $\bullet$ {\bf Solving mean-field control through master Bellman equation and symmetric neural networks.} 
 
\vspace{2mm}

\noindent \cite{GLPW21} solves the master Bellman equation arising from dynamic programming principle applied to mean-field control problems (see \cite{PW17}). The paper approximates the value function evaluated on the empirical measure stemming from particles simulation of a training forward process. The symmetry between iid particles is enforced by optimizing over exchangeable high-dimensional neural networks, invariant by permutation of their inputs. The companion paper \cite{GPW21b} provides a rate for the particle method convergence. 

\vspace{2mm}

\noindent $\bullet$ {\bf Reinforcement Learning for mean-field control and mean-field games \cite{carmona2019modelfree, anahtarc2019fitted,angiuli2020unified,gu2020qlearning,guo2020general}.} 

\vspace{1mm}

\noindent
Some works focus on similar problems but with unknown dynamics. Thus they rely on trajectories sampled from a simulator and reinforcement learnings-- especially Q-learning-- to estimate the state action value function and optimal control without a model. The idea is to optimize a neural network by relying on a memory of past state action transitions used to train the network in order for it to verify the Bellman equation on samples from memory replay.

\vspace{2mm}

\noindent $\bullet$ {\bf Machine learning framework for solving high-dimensional mean field game and mean field control problems
\cite{Ruthotto9183}} 

\vspace{1mm}

\noindent
This paper focuses on potential mean field games, in which the cost functions depending on the law can be written as the linear functional
derivative of a function with respect to a measure. A Lagrangian method with Deep Galerkin type penalization is used. In this case the potential is approached by a neural network and solving mean-field games amounts to solve an unconstrained optimization problem.

\vspace{2mm}

\noindent $\bullet$ {\bf Deep quantum neural networks  \cite{sak20} }
\vspace{1mm}

\noindent
We briefly mention this work studying the use of deep quantum neural networks which exploit the quantum superposition properties by replacing bits by ``qubits". Promising results are obtained when using these networks for regression in financial contexts such as implied volatility estimation. Future works may study the application of such neural networks to control problems and PDEs.

\vspace{2mm}

\noindent $\bullet$ {\bf Path signature for path-dependent PDE \cite{siskaetal20}} 

\vspace{1mm}

\noindent
This work extends previously developed methods for solving state-dependent PDEs to the linear path-dependent setting coming for instance from the pricing and hedging of path-dependent options. A path-dependent Feynman-Kac representation is numerically computed through a global minimization over neural networks. The authors show that using LSTM networks taking the forward process' path signatures (coming from the rough paths literature) as input yields better results than taking the discretized path as input of a feedforward network.

\printbibliography

\end{document}